\documentclass[12pt]{article}
\usepackage{amsmath}
\usepackage{amssymb}
\usepackage{amscd}
\usepackage{amsfonts}

\newtheorem{theorem}{Theorem}

\newtheorem{corollary}{Corollary}

\newtheorem{example}{Example}
\newtheorem{lemma}{Lemma}
\newtheorem{proposition}{Proposition}
\newtheorem{remark}{Remark}
\numberwithin{equation}{section}
\numberwithin{lemma}{section}
\numberwithin{theorem}{section}
\numberwithin{remark}{section}
\numberwithin{corollary}{section}
\numberwithin{proposition}{section}
\numberwithin{definition}{section}
\numberwithin{example}{section}
\numberwithin{table}{section}
\allowdisplaybreaks
\def\func#1{\mathop{\rm #1}}

\begin{document}

\title{Approximation of deterministic and stochastic Navier-Stokes equations
in vorticity-velocity formulation}
\author{G.N. Milstein\thanks{%
Ural Federal University, Lenin Str.~51, 620083 Ekaterinburg, Russia;} \and %
M.V. Tretyakov\thanks{%
School of Mathematical Sciences, University of Nottingham, University Park,
Nottingham, NG7 2RD, UK, email: Michael.Tretyakov@nottingham.ac.uk}}
\maketitle

\begin{abstract}
We consider a time discretization of incompressible Navier-Stokes equations
with spatial periodic boundary conditions in the vorticity-velocity
formulation. The approximation is based on freezing the velocity on time
subintervals resulting in linear parabolic equations for vorticity.
Probabilistic representations for solutions of these linear equations are
given. At each time step, the velocity is expressed via vorticity using a
formula corresponding to the Biot--Savart-type law. We show that the
approximation is divergent free and of first order. The results are extended
to two-dimensional stochastic Navier-Stokes equations with additive noise,
where, in particular, we prove the first mean-square convergence order of
the vorticity approximation.


\end{abstract}

\section{Introduction}

Navier-Stokes equations (NSE), both deterministic and stochastic, are
important for a number of applications and, consequently, development and
analysis of numerical methods for simulation of NSE are of significant
interest. The theory and applications of deterministic NSE can be found,
e.g. in \cite{CF,DJT,MB,FMRT01,RT,T} and for stochastic NSE -- e.g. in \cite%
{Flandoli,MatPhD,RozNS05}. The literature on numerics for deterministic NSE
is extensive \cite{FEM,Peyret,W}\ (see also references therein) while the
literature on numerics for stochastic NSE is still rather sparse, let us
mention \cite{CP,BCP10,BBM14,Dorsek,MTsns}.

In this paper we consider incompressible NSE in the vorticity-velocity
formulation with periodic boundary conditions (see, e.g. \cite{MB} for the
deterministic case and \cite{HaMa06} for the stochastic one). In the
deterministic case we deal with both two-dimensional and three-dimensional
NSE, in the stochastic case we are interested in two-dimensional NSE with
additive noise. We study their time discretization which is based on
freezing the velocity at every time step. Consequently, at every step we
just need to solve a system of linear parabolic PDEs. To compute the
velocity, we express it via the vorticity field, i.e. we derive a periodic
version of Biot-Savart's law (see e.g. \cite[p. 50 and p. 71]{MB}). The
constructed approximations of both vorticity and velocity are divergent
free. We prove convergence theorems for the suggested approximation. The
second part of the paper deals with NSE with additive noise.

The paper is organised as follows. We introduce function spaces required,
recall the Helmholtz-Hodge decomposition, and derive the periodic version of
Biot-Savart's law in Section~\ref{Secpre}. In the deterministic case we
suggest to use probabilistic representations together with ideas of
weak-sense numerical integration of SDEs for solving the system of linear
parabolic PDEs at every step of the time discretization. In Section~\ref%
{SecPrb},\ we present probabilistic representations appropriate for this
task which are based on \cite{M} (see also \cite{MSS,MT1}). The numerical
method is proposed and analysed in Section~\ref{sec:approx}, where in
particular its first order convergence in $L^{2}$-norm is proved. Then the
ideas of Section~\ref{sec:approx} are transferred over to the case of NSE
with additive noise in Section~\ref{sec:sns}, including a proof of
first-order mean-square convergence of the time-discretization of the
stochastic NSE in the vorticity-velocity formulation. Ideas used in the
proof are of potential interest for convergence analysis of numerical
methods for a wider class of semilinear SPDEs.

\section{Preliminaries\label{Secpre}}

In the first part of this paper (Sections~\ref{SecPrb} and~\ref{sec:approx}%
), we consider the two- and three-dimensional deterministic incompressible
NSE for velocity $v$ and pressure $q$ with space periodic conditions:
\begin{eqnarray}
\frac{\partial v}{\partial s}+(v,\nabla )v+\nabla q-\frac{\sigma ^{2}}{2}%
\Delta v=F,  \label{NS1} \\
\mathop{\rm div}v=0.  \label{NS2}
\end{eqnarray}%
In (\ref{NS1})-(\ref{NS2}) we have $0<s\leq T,$\ $x\in \mathbf{R}^{n},$\ $%
v\in \mathbf{R}^{n},$\ $F$ $\in \mathbf{R}^{n},$ $n=2,3,$\ $q$ is a scalar.
The velocity vector $v=(v^{1},\ldots ,v^{n})^{\top }$ satisfies initial
conditions%
\begin{equation}
v(0,x)=\varphi (x)  \label{NS3}
\end{equation}%
and spatial periodic conditions%
\begin{equation}
v(s,x+Le_{i})=v(s,x),\ i=1,\ldots ,n,\ 0\leq s\leq T.  \label{NS4}
\end{equation}%
Here $\mathop{\rm div}\varphi =0,$ $\{e_{i}\}$ is the canonical basis in $%
\mathbf{R}^{n},$ and $L>0$ is the period. For simplicity in writing, the
periods in all the directions are taken to be equal. The function $F=F(s,x)$
and pressure $q=q(s,x)$ are assumed to be spatial periodic as well.

In what follows we will consider the deterministic NSE with negative
direction of time which is convenient for probabilistic representations
considered in Section~\ref{SecPrb}. By an appropriate change of the time
variable and functions, the NSE (\ref{NS1})-(\ref{NS2}) with positive
direction of time can be rewritten in the form:
\begin{eqnarray}
\frac{\partial u}{\partial t}+\frac{\sigma ^{2}}{2}\Delta u-(u,\nabla
)u-\nabla p+f=0,  \label{ns1p} \\
\mathop{\rm div}u=0,  \label{ns2p}
\end{eqnarray}%
where $0\leq t<T,$\ $x\in \mathbf{R}^{n},$\ $u\in \mathbf{R}^{n},$\ $f$ $\in
\mathbf{R}^{n},$ $n=2,3,$\ the pressure $p$ is a scalar. The velocity vector
$u=(u^{1},\ldots ,u^{n})^{\top }$ satisfies the terminal condition%
\begin{equation}
u(T,x)=\varphi (x)  \label{ns3p}
\end{equation}%
and spatial periodic conditions%
\begin{equation}
u(t,x+Le_{i})=u(t,x),\ i=1,\ldots ,n,\ 0\leq t\leq T.  \label{ns4p}
\end{equation}%
Throughout the paper we will assume that this problem has a unique,
sufficiently smooth classical solution. In the two-dimensional case ($n=2)$
the corresponding theory is available, e.g. in \cite{FMRT01,MB}. In the
remaining part of this preliminary section we recall the required function
spaces \cite{DJT,RT,T} and write the NSE in vorticity formulation.

\subsection{Function spaces and the Helmholtz-Hodge decomposition}

Let $\{e_{i}\}$ be the canonical basis in $\mathbf{R}^{n}.$ We shall
consider spatial periodic $n$-vector functions $u(x)=(u^{1}(x),\ldots
,u^{n}(x))^{\top }$ in $\mathbf{R}^{n}:$ $u(x+Le_{i})=u(x),\ i=1,\ldots ,n,$
where $L>0$ is the period in $i$th direction. Denote by $Q=(0,L]^{n}$ the
cube of the period. We denote by $\mathbf{L}^{2}(Q)$ the Hilbert space of
functions on $Q$ with the scalar product and the norm%
\begin{equation*}
(u,v)=\int_{Q}\sum_{i=1}^{n}u^{i}(x)v^{i}(x)dx,\ \Vert u\Vert =(u,u)^{1/2}.
\end{equation*}%
\ \

We keep the notation $|\cdot|$ for the absolute value of numbers and for the
length of $n$-dimensional vectors, for example,%
\begin{equation*}
|u(x)|=[(u^{1}(x))^{2}+\cdots+(u^{n}(x))^{2}]^{1/2}.
\end{equation*}

We denote by $\mathbf{H}_{p}^{m}(Q),\ m=0,1,\ldots,$ the Sobolev space of
functions which are in $\mathbf{L}^{2}(Q),$ together with all their
derivatives of order less than or equal to $m,$ and which are periodic
functions with the period $Q.$ The space $\mathbf{H}_{p}^{m}(Q)$ is a
Hilbert space with the scalar product and the norm%
\begin{equation*}
(u,v)_{m}=\int_{Q}\sum_{i=1}^{n}\sum_{[\alpha^{i}]\leq
m}D^{\alpha^{i}}u^{i}(x)D^{\alpha^{i}}v^{i}(x)dx,\ \Vert
u\Vert_{m}=[(u,u)_{m}]^{1/2},
\end{equation*}
where $\alpha^{i}=(\alpha_{1}^{i},\ldots,\alpha_{n}^{i}),\
\alpha_{j}^{i}\in\{0,\ldots,m\},\
[\alpha^{i}]=\alpha_{1}^{i}+\cdots+\alpha_{n}^{i},$ and%
\begin{equation*}
D^{\alpha^{i}}=D_{1}^{\alpha_{1}^{i}}\cdots D_{n}^{\alpha_{n}^{i}}=\frac{%
\partial^{\lbrack\alpha^{i}]}}{\partial(x^{1})^{\alpha_{1}^{i}}\cdots%
\partial(x^{n})^{\alpha_{n}^{i}}}\ ,\ i=1,\ldots,n.
\end{equation*}
Note that $\mathbf{H}_{p}^{0}(Q)=\mathbf{L}^{2}(Q).$

Introduce the Hilbert subspaces of $\mathbf{H}_{p}^{m}(Q):$%
\begin{eqnarray*}
\mathbf{V}_{p}^{m}& = \{v:\ v\in \mathbf{H}_{p}^{m}(Q),\ \mathop{\rm div}%
v=0\},\ m>0, \\
\mathbf{V}_{p}^{0}& = \text{the closure of } \mathbf{V}_{p}^{m},\ m>0 \text{
in }\mathbf{L}^{2}(Q).
\end{eqnarray*}

Denote by $P$ the orthogonal projection in $\mathbf{H}_{p}^{m}(Q)$ onto $%
\mathbf{V}_{p}^{m}$ (we omit $m$ in the notation $P$ here). The operator $P$
is often called the Leray projection. Due to the Helmholtz-Hodge
decomposition, any function $u\in\mathbf{H}_{p}^{m}(Q)$ can be represented
as
\begin{equation*}
u=Pu+\nabla g,\ \mathop{\rm div}Pu=0,
\end{equation*}
where $g=g(x)$ is a scalar $Q$-periodic function such that $\nabla g\in%
\mathbf{H}_{p}^{m}(Q).$ It is natural to introduce the notation $P^{\bot
}u:=\nabla g$ and hence write
\begin{equation*}
u=Pu+P^{\bot}u
\end{equation*}
with
\begin{equation*}
P^{\bot}u\in(\mathbf{V}_{p}^{m})^{\bot}=\{v:\ v\in\mathbf{H}_{p}^{m}(Q),\
v=\nabla g\}.
\end{equation*}

Let%
\begin{eqnarray}
u(x)=\sum_{\mathbf{n}\in \mathbf{Z}^{n}}u_{\mathbf{n}}e^{i(2\pi /L)(\mathbf{n%
},x)},\ g(x)=\sum_{\mathbf{n}\in \mathbf{Z}^{n}}g_{\mathbf{n}}e^{i(2\pi /L)(%
\mathbf{n},x)},\ g_{\mathbf{0}}=0,  \label{N00} \\
Pu(x)=\sum_{\mathbf{n}\in \mathbf{Z}^{n}}(Pu)_{\mathbf{n}}e^{i(2\pi /L)(%
\mathbf{n},x)},\ P^{\bot }u(x)=\nabla g(x)=\sum_{\mathbf{n}\in \mathbf{Z}%
^{n}}(P^{\bot }u)_{\mathbf{n}}e^{i(2\pi /L)(\mathbf{n},x)}  \notag
\end{eqnarray}%
be the Fourier expansions of $u,$\ $g,$\ $Pu,$ and $P^{\bot }u=\nabla g.$
Here $u_{\mathbf{n}},$\ $(Pu)_{\mathbf{n}},\ $and $(P^{\bot }u)_{\mathbf{n}%
}=(\nabla g)_{\mathbf{n}}$ are $n$-dimensional vectors and $g_{\mathbf{n}}$
are scalars. We note that $g_{\mathbf{0}}$ can be any real number but for
definiteness we set $g_{\mathbf{0}}=0$ without loss of generality \cite%
{FMRT01}. The coefficients $(Pu)_{\mathbf{n}},\ (P^{\bot }u)_{\mathbf{n}}$,
and $g_{\mathbf{n}}$ can be easily expressed in terms of $u_{\mathbf{n}}:$
\begin{eqnarray}
(Pu)_{\mathbf{n}}& =u_{\mathbf{n}}-\frac{u_{\mathbf{n}}^{\top }\mathbf{n}}{|%
\mathbf{n}|^{2}}\mathbf{n,\ }(P^{\bot }u)_{\mathbf{n}}=i\frac{2\pi }{L}g_{%
\mathbf{n}}\mathbf{n=}\frac{u_{\mathbf{n}}^{\top }\mathbf{n}}{|\mathbf{n}%
|^{2}}\mathbf{n,\ }  \label{N01} \\
g_{\mathbf{n}}& =-i\frac{L}{2\pi }\frac{u_{\mathbf{n}}^{\top }\mathbf{n}}{|%
\mathbf{n}|^{2}},\ \mathbf{n\neq 0,\ }g_{\mathbf{0}}=0.  \notag
\end{eqnarray}

We have%
\begin{equation*}
\nabla e^{i(2\pi /L)(\mathbf{n},x)}=\mathbf{n}e^{i(2\pi /L)(\mathbf{n}%
,x)}\cdot i\frac{2\pi }{L},
\end{equation*}%
hence $u_{\mathbf{n}}e^{i(2\pi /L)(\mathbf{n},x)}\in \mathbf{V}_{p}^{m}$ if
and only if $(u_{\mathbf{n}},\mathbf{n)}=0.$ We obtain from here that the
orthogonal basis of the subspace $(\mathbf{V}_{p}^{m})^{\bot }$ consists of $%
\mathbf{n}e^{i(2\pi /L)(\mathbf{n},x)},\ \mathbf{n}\in \mathbf{Z}^{n},\
\mathbf{n\neq 0}$; and an orthogonal basis of $\mathbf{V}_{p}^{m}$ consists
of $_{k}u_{\mathbf{n}}e^{i(2\pi /L)(\mathbf{n},x)},$\ $k=1,\ldots ,n-1,\
\mathbf{n}\in \mathbf{Z}^{n}$, where under $\mathbf{n\neq 0}$ the vectors $%
_{k}u_{\mathbf{n}}$ are orthogonal to $\mathbf{n:}$\textbf{\ }$\mathbf{(}%
_{k}u_{\mathbf{n}},\mathbf{n)}=0,\ k=1,\ldots ,n-1,$ and they are orthogonal
among themselves: $\mathbf{(}_{k}u_{\mathbf{n}},\ _{m}u_{\mathbf{n}}\mathbf{)%
}=0,$\ $k,m=1,\ldots ,n-1,$\ $m\neq k,$ and finally, for $\mathbf{n=0,}$ the
vectors $_{k}u_{\mathbf{0}},\ k=1,\ldots ,n,$ are orthogonal. In particular,
in the two-dimensional case $(n=2),$ these bases are, correspondingly (for $%
\mathbf{n}\neq 0$):
\begin{equation}
\left[
\begin{array}{c}
n_{1} \\
n_{2}%
\end{array}%
\right] e^{i(2\pi /L)(\mathbf{n},x)} \text{ \ \ and \ }\ \left[
\begin{array}{c}
-n_{2} \\
n_{1}%
\end{array}%
\right] e^{i(2\pi /L)(\mathbf{n},x)},\ \mathbf{n}=(n_{1},n_{2})^{\top }.
\label{N02}
\end{equation}

We shall consider the case of zero space average (see e.g. \cite{FMRT01}),
i.e. when
\begin{equation}
\int_{Q}u(x)=0.  \label{av}
\end{equation}%
In this case the Fourier series expansion for $u(x)$ does not contain the
constant term and $\sum_{\mathbf{n}\in \mathbf{Z}^{n}}$ in (\ref{N00}) can
be replaced by $\sum_{\mathbf{n}\in \mathbf{Z}^{n},\mathbf{n}\neq 0},$ which
in what follows we will write as simply $\sum .$

We recall Parseval's identity
\begin{equation}
\Vert u\Vert ^{2}=\int_{Q}|u(x)|^{2}dx=L^{n}\sum |u_{\mathbf{n}}|^{2},\ \
n=2,3.  \label{eq:pars}
\end{equation}%
We also note the following two relationships. Since the vector field $u=u(x)$
is real valued, we have
\begin{equation*}
u_{-\mathbf{n}}=\bar{u}_{\mathbf{n}},\ \ \mathbf{n}\in \mathbf{Z}^{n},\
\mathbf{n}\neq 0,
\end{equation*}%
where $\bar{u}_{\mathbf{n}}$ denotes the complex conjugate of $u_{\mathbf{n}%
}.$ The divergence-free condition reads
\begin{equation*}
u_{\mathbf{n}}^{\top }\mathbf{n}=(u_{\mathbf{n}},\mathbf{n}%
)=\sum_{k=1}^{n}u_{\mathbf{n}}^{k}\mathbf{n}^{k}=0,\ \ n=2,3.
\end{equation*}

We will need the following estimate for the tri-linear form (see \cite[p.
50, eq. (6.10)]{CF} or \cite[p. 12, eq. (2.29)]{T}):
\begin{equation}
|((v,\nabla )u,g)|\leq K\Vert v\Vert _{m_{1}}\Vert u\Vert _{m_{2}+1}\Vert
g\Vert _{m_{3}},  \label{A2}
\end{equation}%
where $K>0$ is a constant, $m_{1},\ m_{2}$ and $m_{3}$ are such that $%
m_{1}+m_{2}+m_{3}\geq n/2$ and $(m_{1},m_{2},m_{3})\neq (0,0,n/2),$ $%
(0,n/2,0),$ $(n/2,0,0),$ and $u,$ $v,$ $g$ are arbitrary functions from the
corresponding spaces. Further, we recall the standard interpolation
inequality for Sobolev spaces (see e.g. \cite[p. 11]{T}):
\begin{equation}
\Vert u\Vert _{m}\leq \Vert u\Vert _{m_{1}}^{1-l}\Vert u\Vert _{m_{2}}^{l},
\label{A4i}
\end{equation}%
where $m=(1-l)m_{1}+lm_{2},\ m_{1},$ $m_{2}\geq 0,$ $l\in (0,1),$ and $u\in
\mathbf{H}_{p}^{\max (m_{1},m_{2})}(Q).$ For any $c>0,$ we get from (\ref%
{A4i}) and Young's inequality:
\begin{equation}
\Vert u\Vert _{m}^{2}\leq \Vert u\Vert _{m_{1}}^{2-2l}\Vert u\Vert
_{m_{2}}^{2l}\leq (1-l)c\Vert u\Vert _{m_{1}}^{2}+lc^{1-\frac{1}{l}}\Vert
u\Vert _{m_{2}}^{2}.  \label{A4}
\end{equation}%
Let us take $m_{1}=1,$ $m_{2}=0$ and $m_{3}=1/2$ in (\ref{A2}), then, using (%
\ref{A4}), we get (see also \cite[p. 1028, eq. (A6)]{HaMa06}) for any $%
c_{1}>0$ and $c_{2}>0$:
\begin{eqnarray}
|((v,\nabla )u,g)| &\leq &K\Vert v\Vert _{1}\Vert u\Vert _{1}\Vert g\Vert
_{1/2}\leq \frac{K^{2}}{4c_{1}}\Vert g\Vert _{1/2}^{2}+c_{1}\Vert v\Vert
_{1}^{2}\Vert u\Vert _{1}^{2}  \label{A6} \\
&\leq &\frac{K^{2}}{4c_{1}}\Vert g\Vert _{1/2}^{2}+c_{1}\Vert v\Vert
_{1}^{2}\Vert u\Vert _{1}^{2}  \notag \\
&\leq &c_{2}\Vert g\Vert _{1}^{2}+\frac{K^{4}}{64c_{1}^{2}c_{2}}\Vert g\Vert
^{2}+c_{1}\Vert v\Vert _{1}^{2}\Vert u\Vert _{1}^{2},  \notag
\end{eqnarray}%
where we used (\ref{A4}) with $m=1/2,$ $m_{1}=1,$ $m_{2}=0$ and $l=1/2$
(fractional $\mathbf{H}_{p}^{m}(Q)$ spaces are defined in the usual way via
the Fourier series expansions, see e.g. \cite[pp. 7-8]{T}).

We also recall (see e.g. \cite[p. 20 \ eq. (4.14) ]{FMRT01}) Poincare's
inequality for functions $u\in \mathbf{H}_{p}^{1}(Q)$ satisfying (\ref{av}):
\begin{equation}
||u||\leq \alpha ||\nabla u||  \label{Poin}
\end{equation}%
for some constant $\alpha >0$ which depends only on the shape of $Q$ and on
the period $L.$ We note that here and in what follows: when $u(x)$ is a
vector, $\nabla u(x)$ means the matrix with elements $\partial
u^{i}/\partial x^{j}$ and $||\nabla u||$ means $L^{2}$-norm of the Frobenius
norm of the matrix $\nabla u(x).$

\subsection{Equation for vorticity}

Introduce the vorticity $\omega :$%
\begin{equation}
\omega =\left[
\begin{array}{c}
\omega ^{1} \\
\omega ^{2} \\
\omega ^{3}%
\end{array}%
\right] :=\mathop{\rm curl}u=\mathop{\rm curl}\left[
\begin{array}{c}
u^{1} \\
u^{2} \\
u^{3}%
\end{array}%
\right] =\left[
\begin{array}{ccc}
i & j & k \\
\frac{\partial }{\partial x^{1}} & \frac{\partial }{\partial x^{2}} & \frac{%
\partial }{\partial x^{3}} \\
u^{1} & u^{2} & u^{3}%
\end{array}%
\right] =\left[
\begin{array}{c}
\frac{\partial u^{3}}{\partial x^{2}}-\frac{\partial u^{2}}{\partial x^{3}}
\\
\frac{\partial u^{1}}{\partial x^{3}}-\frac{\partial u^{3}}{\partial x^{1}}
\\
\frac{\partial u^{2}}{\partial x^{1}}-\frac{\partial u^{1}}{\partial x^{2}}%
\end{array}%
\right] .  \label{vor1}
\end{equation}%
We note that (\ref{vor1}) implies $\mathop{\rm div}\omega =0.$

Taking the $\mathop{\rm curl}$\ of equation (\ref{ns1p}) gives the evolution
equation for the vorticity $\omega =\mathop{\rm curl}u:$
\begin{equation}
\frac{\partial \omega }{\partial t}-(u,\nabla )\omega +(\omega ,\nabla )u+%
\frac{\sigma ^{2}}{2}\Delta \omega +g=0,  \label{V1}
\end{equation}%
where $g=\mathop{\rm curl}f.$ From (\ref{ns3p})-(\ref{ns4p}), we get
\begin{equation}
\omega (T,x)=\mathop{\rm curl}\varphi (x):=\phi (x)
\end{equation}%
and spatial periodic conditions%
\begin{equation}
\omega (t,x+Le_{i})=\omega (t,x),\ i=1,\ldots ,n,\ 0\leq t\leq T.
\end{equation}

Analogously to (\ref{N00}), we write the Fourier expansion for $\omega :$
\begin{equation}
\omega (t,x)=\sum \omega _{\mathbf{n}}(t)e^{i(2\pi /L)(\mathbf{n},x)},
\label{vor2}
\end{equation}%
where
\begin{eqnarray*}
\omega _{\mathbf{n}}(t) &=&\frac{1}{L^{n}}\left( \omega (t,\cdot ),e^{i(2\pi
/L)(\mathbf{n},\cdot )}\right) \\
&=&\frac{1}{L^{n}}\int_{Q}\omega (t,x)e^{-i(2\pi /L)(\mathbf{n},x)}dx,\ \
\mathbf{n\in Z}^{n},\ \mathbf{n}\neq 0.
\end{eqnarray*}

Substituting the Fourier expansions for $\omega $ and $u$ in (\ref{vor1}),
we obtain
\begin{equation}
\omega =\left[
\begin{array}{c}
\sum \omega _{\mathbf{n}}^{1}(t)e^{i(2\pi /L)(\mathbf{n},x)} \\
\sum \omega _{\mathbf{n}}^{2}(t)e^{i(2\pi /L)(\mathbf{n},x)} \\
\sum \omega _{\mathbf{n}}^{3}(t)e^{i(2\pi /L)(\mathbf{n},x)}%
\end{array}%
\right] =i\frac{2\pi }{L}\left[
\begin{array}{c}
\sum \left( u_{\mathbf{n}}^{3}(t)\mathbf{n}^{2}-u_{\mathbf{n}}^{2}(t)\mathbf{%
n}^{3}\right) e^{i(2\pi /L)(\mathbf{n},x)} \\
\sum \left( u_{\mathbf{n}}^{1}(t)\mathbf{n}^{3}-u_{\mathbf{n}}^{3}(t)\mathbf{%
n}^{1}\right) e^{i(2\pi /L)(\mathbf{n},x)} \\
\sum \left( u_{\mathbf{n}}^{2}(t)\mathbf{n}^{1}-u_{\mathbf{n}}^{1}(t)\mathbf{%
n}^{2}\right) e^{i(2\pi /L)(\mathbf{n},x)}%
\end{array}%
\right] .  \label{vor3}
\end{equation}%
The equality (\ref{vor3}) gives for any $\mathbf{n}\neq 0$ the equations
with respect to $u_{\mathbf{n}}^{k},$ $k=1,2,3:$
\begin{eqnarray}
\mathbf{n}^{2}u_{\mathbf{n}}^{3}-\mathbf{n}^{3}u_{\mathbf{n}}^{2} &=&-\frac{%
iL}{2\pi }\omega _{\mathbf{n}}^{1}  \label{vor4} \\
\mathbf{n}^{3}u_{\mathbf{n}}^{1}-\mathbf{n}^{1}u_{\mathbf{n}}^{3} &=&-\frac{%
iL}{2\pi }\omega _{\mathbf{n}}^{2}  \notag \\
\mathbf{n}^{1}u_{\mathbf{n}}^{2}-\mathbf{n}^{2}u_{\mathbf{n}}^{1} &=&-\frac{%
iL}{2\pi }\omega _{\mathbf{n}}^{3}\ .  \notag
\end{eqnarray}%
Thanks to $\mathop{\rm div}u=0$ and $\mathop{\rm div}\omega =0,$ we also
have for any $\mathbf{n:}$
\begin{eqnarray}
\mathbf{n}^{1}u_{\mathbf{n}}^{1}+\mathbf{n}^{2}u_{\mathbf{n}}^{2}+\mathbf{n}%
^{3}u_{\mathbf{n}}^{3} &=&0  \label{vor5} \\
\mathbf{n}_{\mathbf{n}}^{1}\omega _{\mathbf{n}}^{1}+\mathbf{n}^{2}\omega _{%
\mathbf{n}}^{2}+\mathbf{n}^{3}\omega _{\mathbf{n}}^{3} &=&0\ .  \label{vor6}
\end{eqnarray}%
Due to the property (\ref{vor6}), the system (\ref{vor4})-(\ref{vor5}) with
respect to $u_{\mathbf{n}}^{k}$ is compatible. It is not difficult to prove
directly that the solution of this system is unique. This also follows from
the two observations that the vector field $u_{\mathbf{n}}e^{i(2\pi /L)(%
\mathbf{n},x)}$ is solenoidal because it is divergence-free and that in the
case of $\omega _{\mathbf{n}}^{k}=0,$ $k=1,2,3,$ the field is irrotational,
i.e. potential. But if a vector field is simultaneously solenoidal and
potential, it is trivial, i.e. $u_{\mathbf{n}}=0.$ Thus, the homogeneous
system corresponding to (\ref{vor4})-(\ref{vor5}) has the trivial solution
only, and hence the solution to the system (\ref{vor4})-(\ref{vor5}) exists
and it is unique.

Our nearest goal consists in solving this system, i.e., in expressing $u$
via $\omega .$ We observe that
\begin{eqnarray}
\mathop{\rm curl}u &=&\omega ,  \label{vor7} \\
\mathop{\rm div}u &=&0.  \label{vor8}
\end{eqnarray}

\begin{proposition}
\label{prop:curl} For a sufficiently smooth $\psi ,$ let $\mathop{\rm div}%
\psi =0.$ Then%
\begin{equation}
\mathop{\rm curl}[-\func{curl}\psi ]=\Delta \psi .  \label{V5}
\end{equation}
\end{proposition}

\noindent \textbf{Proof}. The first component of the vector $%
\mathop{\rm
curl}[-\func{curl}\psi ]$ is equal to%
\begin{eqnarray*}
&&-\frac{\partial ^{2}\psi ^{2}}{\partial x^{1}\partial x^{2}}+\frac{%
\partial ^{2}\psi ^{1}}{(\partial x^{2})^{2}}+\frac{\partial ^{2}\psi ^{1}}{%
(\partial x^{3})^{2}}-\frac{\partial ^{2}\psi ^{3}}{\partial x^{1}\partial
x^{3}} \\
&=&\frac{\partial }{\partial x^{1}}(-\frac{\partial \psi ^{2}}{\partial x^{2}%
}-\frac{\partial \psi ^{3}}{\partial x^{3}})+\frac{\partial ^{2}\psi ^{1}}{%
(\partial x^{2})^{2}}+\frac{\partial ^{2}\psi ^{1}}{(\partial x^{3})^{2}}.
\end{eqnarray*}%
Because of the condition $\mathop{\rm div}\psi =0,$ this component is equal
to $\Delta \psi ^{1}.$ Analogously, the second and third components are
equal to $\Delta \psi ^{2}$ and $\Delta \psi ^{3},$ correspondingly. The
proposition is proved.\bigskip

Let us look for $u$\ in the form%
\begin{equation}
u=-\mathop{\rm curl}\psi ,  \label{V7}
\end{equation}%
where $\mathop{\rm div}\psi =0.$ Due to (\ref{vor7})-(\ref{V5}), we now have
to solve%
\begin{equation}
\Delta \psi =\omega ,\ \omega =\sum \omega _{\mathbf{n}}e^{i(2\pi /L)(%
\mathbf{n},x)}.  \label{V8}
\end{equation}%
Equation (\ref{V8}) is solvable (uniquely, if we assume $\psi _{\mathbf{0}%
}=0 $):%
\begin{equation*}
\psi =\sum \psi _{\mathbf{n}}e^{i(2\pi /L)(\mathbf{n},x)},\ \psi _{\mathbf{n}%
}=-\frac{\omega _{\mathbf{n}}L^{2}}{4\pi ^{2}|\mathbf{n}|^{2}},
\end{equation*}%
i.e.,
\begin{equation*}
\psi _{\mathbf{n}}^{j}=-\frac{\omega _{\mathbf{n}}^{j}L^{2}}{4\pi ^{2}|%
\mathbf{n}|^{2}},\ j=1,2,3.
\end{equation*}%
Hence, using (\ref{V7}), we have%
\begin{equation}
u=-\mathop{\rm curl}\psi =\frac{Li}{2\pi }\left[
\begin{array}{c}
\sum \frac{1}{|\mathbf{n}|^{2}}e^{i(2\pi /L)(\mathbf{n},x)}(\omega _{\mathbf{%
n}}^{3}\mathbf{n}^{2}-\omega _{\mathbf{n}}^{2}\mathbf{n}^{3}) \\
\sum \frac{1}{|\mathbf{n}|^{2}}e^{i(2\pi /L)(\mathbf{n},x)}(\omega _{\mathbf{%
n}}^{1}\mathbf{n}^{3}-\omega _{\mathbf{n}}^{3}\mathbf{n}^{1}) \\
\sum \frac{1}{|\mathbf{n}|^{2}}e^{i(2\pi /L)(\mathbf{n},x)}(\omega _{\mathbf{%
n}}^{2}\mathbf{n}^{1}-\omega _{\mathbf{n}}^{1}\mathbf{n}^{2})%
\end{array}%
\right] :=U\omega ,  \label{V9}
\end{equation}%
where $U$ is a linear operator. It is not difficult to verify that the
equality $\mathop{\rm div}u=0$ (see (\ref{vor8})) holds for $u$\ from (\ref%
{V9}) under arbitrary $\omega $. However, the equality (\ref{vor7}) is not
fulfilled by $u$\ from (\ref{V9}) for arbitrary $\omega $. But the
considered $\omega $ is not arbitrary, it is divergent free and we show
below that for a divergent-free $\omega $ the equality (\ref{vor7}) is
satisfied by $u$ from (\ref{V9}).

For $u$\ from (\ref{V9}), we get%
\begin{equation*}
\mathop{\rm curl}u=-\left[
\begin{array}{c}
\sum \frac{1}{|\mathbf{n}|^{2}}e^{i(2\pi /L)(\mathbf{n},x)}[-\omega _{%
\mathbf{n}}^{1}((\mathbf{n}^{2})^{2}+(\mathbf{n}^{3})^{2})+\mathbf{n}^{1}(%
\mathbf{n}^{2}\omega _{\mathbf{n}}^{2}+\mathbf{n}^{3}\omega _{\mathbf{n}%
}^{3})] \\
\sum \frac{1}{|\mathbf{n}|^{2}}e^{i(2\pi /L)(\mathbf{n},x)}[-\omega _{%
\mathbf{n}}^{2}((\mathbf{n}^{1})^{2}+(\mathbf{n}^{3})^{2})+\mathbf{n}^{2}(%
\mathbf{n}^{1}\omega _{\mathbf{n}}^{1}+\mathbf{n}^{3}\omega _{\mathbf{n}%
}^{3})] \\
\sum \frac{1}{|\mathbf{n}|^{2}}e^{i(2\pi /L)(\mathbf{n},x)}[-\omega _{%
\mathbf{n}}^{3}((\mathbf{n}^{1})^{2}+(\mathbf{n}^{2})^{2})+\mathbf{n}^{3}(%
\mathbf{n}^{1}\omega _{\mathbf{n}}^{1}+\mathbf{n}^{2}\omega _{\mathbf{n}%
}^{2})]%
\end{array}%
\right] .
\end{equation*}%
Recall that (\ref{vor6}) holds for a divergent-free $\omega .$ Then, using (%
\ref{vor6}), we have%
\begin{eqnarray*}
\mathbf{n}^{1}(\mathbf{n}^{2}\omega _{\mathbf{n}}^{2}+\mathbf{n}^{3}\omega _{%
\mathbf{n}}^{3}) &=&-(\mathbf{n}^{1})^{2}\omega _{\mathbf{n}}^{1} \\
\mathbf{n}^{2}(\mathbf{n}^{1}\omega _{\mathbf{n}}^{1}+\mathbf{n}^{3}\omega _{%
\mathbf{n}}^{3}) &=&-(\mathbf{n}^{2})^{2}\omega _{\mathbf{n}}^{2} \\
\mathbf{n}^{3}(\mathbf{n}^{1}\omega _{\mathbf{n}}^{1}+\mathbf{n}^{2}\omega _{%
\mathbf{n}}^{2}) &=&-(\mathbf{n}^{3})^{2}\omega _{\mathbf{n}}^{3}
\end{eqnarray*}%
and therefore%
\begin{equation*}
\mathop{\rm curl}u=\left[
\begin{array}{c}
\sum \frac{1}{|\mathbf{n}|^{2}}e^{i(2\pi /L)(\mathbf{n},x)}\omega _{\mathbf{n%
}}^{1}|\mathbf{n}|^{2} \\
\sum \frac{1}{|\mathbf{n}|^{2}}e^{i(2\pi /L)(\mathbf{n},x)}\omega _{\mathbf{n%
}}^{2}|\mathbf{n}|^{2} \\
\sum \frac{1}{|\mathbf{n}|^{2}}e^{i(2\pi /L)(\mathbf{n},x)}\omega _{\mathbf{n%
}}^{3}|\mathbf{n}|^{2}%
\end{array}%
\right] ,
\end{equation*}%
i.e., $\mathop{\rm curl}u=\omega .$ Thus, the following theorem is proved.

\begin{theorem}
\label{thm:nsvor}The velocity field $u$\ (with $u_{\mathbf{0}}=0$) is
determined explicitly through vorticity field $\omega =\mathop{\rm curl}v$
by formula $(\ref{V9})$. The closed-form equation for $\omega $ is given by%
\begin{equation}
\frac{\partial \omega }{\partial t}-(U\omega ,\nabla )\omega +(\omega
,\nabla )U\omega +\frac{\sigma ^{2}}{2}\Delta \omega +g=0,  \label{V11}
\end{equation}%
where $U$ is from $(\ref{V9})$.
\end{theorem}

This theorem is related to analogous results for periodic $2D$ flows (see
\cite[p. 50]{MB}) and for flows in the whole space (see \cite[ p. 71]{MB}).

Let us mention the corresponding formulas in the 2D case, for which we have%
\begin{equation*}
u(x)=\left[
\begin{array}{c}
u^{1}(x^{1},x^{2}) \\
u^{2}(x^{1},x^{2}) \\
0%
\end{array}%
\right] =\sum \left[
\begin{array}{c}
u_{\mathbf{n}}^{1} \\
u_{\mathbf{n}}^{2} \\
0%
\end{array}%
\right] e^{i(2\pi /L)(\mathbf{n}^{1}x^{1}+\mathbf{n}^{2}x^{2})},
\end{equation*}%
i.e., $u^{1}(x)$\ and $u^{2}(x)$ are independent of $x^{3}$ and $u^{3}=0.$
Hence%
\begin{equation*}
\omega =\left[
\begin{array}{c}
\omega ^{1} \\
\omega ^{2} \\
\omega ^{3}%
\end{array}%
\right] =\left[
\begin{array}{c}
0 \\
0 \\
\frac{\partial u^{2}}{\partial x^{1}}-\frac{\partial u^{1}}{\partial x^{2}}%
\end{array}%
\right] .
\end{equation*}%
We shall denote the scalar $\omega ^{3}(x)$ as $\omega (x)=\frac{\partial
u^{2}}{\partial x^{1}}(x)-\frac{\partial u^{1}}{\partial x^{2}}(x)$ and the
two dimensional vector $(u^{1}(x^{1},x^{2}),u^{2}(x^{1},x^{2}))^{\intercal }$
as $u(x).$ This does not lead to confusion. We have $\mathop{\rm div}u(x)=%
\frac{\partial u^{1}}{\partial x^{1}}(x)+\frac{\partial u^{2}}{\partial x^{2}%
}(x)=0 $, i.e., $u_{\mathbf{n}}^{1}\mathbf{n}^{1}+u_{\mathbf{n}}^{2}\mathbf{n%
}^{2}=0 $ for any $\mathbf{n}=(\mathbf{n}^{1},\mathbf{n}^{2}),$ and%
\begin{equation*}
\omega (x)=\sum \omega _{\mathbf{n}}e^{i(2\pi /L)(\mathbf{n}^{1}x^{1}+%
\mathbf{n}^{2}x^{2})},
\end{equation*}%
where%
\begin{equation*}
\omega _{\mathbf{n}}=i\frac{2\pi }{L}(u_{\mathbf{n}}^{2}\mathbf{n}^{1}-u_{%
\mathbf{n}}^{1}\mathbf{n}^{2}).
\end{equation*}%
Due to (\ref{V9}), $u$ is expressed through $\omega :$%
\begin{equation}
u(x)=\frac{Li}{2\pi }\left[
\begin{array}{c}
\sum \frac{1}{|\mathbf{n}|^{2}}e^{i(2\pi /L)(\mathbf{n},x)}\omega _{\mathbf{n%
}}\mathbf{n}^{2} \\
-\sum \frac{1}{|\mathbf{n}|^{2}}e^{i(2\pi /L)(\mathbf{n},x)}\omega _{\mathbf{%
n}}\mathbf{n}^{1}%
\end{array}%
\right] :=U\omega .  \label{Mnew1}
\end{equation}%
Clearly, $(\omega ,\nabla )u=0$ in the $2D$ case. Hence (\ref{V1}) takes the
form%
\begin{equation}
\frac{\partial \omega }{\partial t}-u^{1}(t,x)\frac{\partial \omega }{%
\partial x^{1}}(t,x)-u^{2}(t,x)\frac{\partial \omega }{\partial x^{2}}(t,x)+%
\frac{\sigma ^{2}}{2}\Delta \omega (t,x)+g(t,x^{1},x^{2})=0.  \label{M6}
\end{equation}

\begin{example}
\label{ex:21} Consider the Stokes equation with negative direction of time:
\begin{equation}
\frac{\partial u}{\partial t}+\frac{\sigma ^{2}}{2}\Delta u-\nabla p+f=0
\label{vor15}
\end{equation}%
with the conditions $(\ref{ns2p})$-$(\ref{ns4p})$. In this case the
vorticity $\omega (t,x)$ satisfies the equation
\begin{equation}
\frac{\partial \omega }{\partial t}+\frac{\sigma ^{2}}{2}\Delta \omega +g=0,
\label{vor16}
\end{equation}%
where $g(t,x)=\sum g_{\mathbf{n}}e^{i(2\pi /L)(\mathbf{n},x)}.$ Substituting
the Fourier expansions for $\omega $ and $g$ in $(\ref{vor16})$, we get%
\begin{equation*}
\sum \omega _{\mathbf{n}}^{\prime }(t)e^{i(2\pi /L)(\mathbf{n},x)}+\frac{%
\sigma ^{2}}{2}\sum \omega _{\mathbf{n}}(t)\left( -\frac{4\pi ^{2}}{L^{2}}%
\right) |\mathbf{n}|^{2}e^{i(2\pi /L)(\mathbf{n},x)}+\sum g_{\mathbf{n}%
}e^{i(2\pi /L)(\mathbf{n},x)}=0,
\end{equation*}%
whence
\begin{equation*}
\frac{d\omega _{\mathbf{n}}^{k}}{dt}-\frac{2\pi ^{2}\sigma ^{2}}{L^{2}}|%
\mathbf{n}|^{2}\omega _{\mathbf{n}}^{k}+g_{\mathbf{n}}^{k}=0,\ \ \omega _{%
\mathbf{n}}^{k}(T)=\phi _{\mathbf{n}}^{k},
\end{equation*}%
where $\phi _{\mathbf{n}}$ are the Fourier coefficients for $\phi :=%
\mathop{\rm curl}\varphi .$ Hence
\begin{equation*}
\omega _{\mathbf{n}}^{k}(t)=\phi _{\mathbf{n}}^{k}\exp \left( \frac{2\pi
^{2}\sigma ^{2}}{L^{2}}|\mathbf{n}|^{2}(t-T)\right) +\int_{t}^{T}\exp \left(
\frac{2\pi ^{2}\sigma ^{2}}{L^{2}}|\mathbf{n}|^{2}(t-s)\right) g_{\mathbf{n}%
}^{k}(s)ds.
\end{equation*}
\end{example}

In future we will need the following estimates. One can obtain from (\ref{V9}%
) that for $m\geq 1$%
\begin{equation}
||u||_{m}=||U\omega ||_{m}\leq K||\omega ||_{m-1}  \label{A3}
\end{equation}%
for some $K>0.$ Further, we note that
\begin{equation}
||\omega ||_{1}^{2}=||\omega ||^{2}+||\nabla \omega ||^{2}  \label{A31}
\end{equation}%
and then by (\ref{Poin})%
\begin{equation}
||\omega ||_{1}^{2}\leq K||\nabla \omega ||^{2}  \label{A32}
\end{equation}%
for some $K>0$. Using (\ref{A6}), (\ref{A3}), (\ref{A31}) and (\ref{A32}),
we get for $\omega ,$ $v$, $g$ from appropriate spaces and arbitrary $%
c_{1}>0 $ and $c_{2}>0:$
\begin{eqnarray}
|((U\omega ,\nabla )v,g)| &\leq &c_{2}\Vert g\Vert _{1}^{2}+\frac{K^{4}}{%
64c_{1}^{2}c_{2}}\Vert g\Vert ^{2}+c_{1}\Vert U\omega \Vert _{1}^{2}\Vert
v\Vert _{1}^{2}  \label{A66} \\
&\leq &c_{2}||\nabla g||^{2}+(c_{2}+\frac{K^{4}}{64c_{1}^{2}c_{2}})\Vert
g\Vert ^{2}+Kc_{1}\Vert \omega \Vert ^{2}\Vert \nabla v\Vert ^{2}  \notag \\
&=&c_{2}||\nabla g||^{2}+K\Vert g\Vert ^{2}+c_{3}\Vert \omega \Vert
^{2}\Vert \nabla v\Vert ^{2},  \notag
\end{eqnarray}%
where in the third line $c_{3}>0$ is an arbitrary constant and $K>0$ is some
constant dependent on $c_{2}$ and $c_{3}$\ (it differs from $K>0$ in the
first and second line but this should not cause any confusion).

\section{Probabilistic representations of solutions to linear systems of
parabolic equations with application to vorticity equations\label{SecPrb}}

In this section we derive probabilistic representations for systems of
parabolic equations based on the approach developed in \cite{M}. They can be
used for constructing probabilistic methods for NSE in vorticity-velocity
formulation (\ref{V1}) (see probabilistic numerical methods for semilinear
PDEs in e.g. \cite{M1,MT1} and for NSE in velocity formulation in e.g. \cite%
{MTns13}). For this purpose, it is useful to have a wide class of such
probabilistic representations, and, in addition to \cite{M}, we also exploit
ideas from \cite{MSS,MT1}. Note that we obtain more general probabilistic
representations than in \cite{BFR}.

\subsection{The basic probabilistic representation}

We consider the following Cauchy problem for system of parabolic equations%
\begin{eqnarray}
\frac{\partial u^{k}}{\partial t}+\frac{1}{2}\sum_{i,j=1}^{n}\left[
\sum_{r=1}^{l}\sigma _{r}\sigma _{r}^{\top }\right] ^{ij}\frac{\partial
^{2}u^{k}}{\partial x^{i}\partial x^{j}}+\sum_{i=1}^{n}\sum_{j=1}^{m}\left(
\sum_{r=1}^{l}[\sigma _{r}]^{i}[\vartheta _{r}]^{jk}\right) \frac{\partial
u^{j}}{\partial x^{i}}  \label{PR1} \\
+\sum_{i=1}^{n}a^{i}\frac{\partial u^{k}}{\partial x^{i}}+[B^{\top
}u]^{k}+f^{k}=0,\ k=1,\ldots ,m,  \notag
\end{eqnarray}%
\begin{equation}
u(T,x)=\varphi (x).  \label{PR2}
\end{equation}

Introduce the system of SDEs%
\begin{eqnarray}
dX &=&a(s,X)ds+\sum_{r=1}^{l}\sigma _{r}(s,X)dw_{r}(s),\ X(t)=x,  \label{PR3}
\\
dY &=&B(s,X)Yds+\sum_{r=1}^{l}\vartheta _{r}(s,X)Ydw_{r}(s),\ Y(t)=y.
\label{PR4}
\end{eqnarray}

In (\ref{PR1})-(\ref{PR4}), $0\leq t\leq s\leq T;\ x$ and $X$ are
column-vectors of dimension $n;\ y$ and $Y$ are column-vectors of dimension $%
m;\ w_{r},\ r=1,\ldots ,l,$ are independent standard Wiener processes; $%
a(s,x)$ and $\sigma _{r}(s,x)$ are column-vectors of dimension $n;\ B(s,x)$
and $\vartheta _{r}(s,x)$ are $m\times m$ - matrices; $u(s,x),\ f(s,x),$ and
$\varphi (x)$ are column-vector of dimension $m$ with components $u^{k},\
f^{k},\ \varphi ^{k},\ k=1,\ldots ,m.$ We assume that there exist a
sufficiently smooth solution of the problem (\ref{PR1})-(\ref{PR2}) and a
unique solution of the problem (\ref{PR3})-(\ref{PR4}).

Introduce the process%
\begin{equation}
\xi _{t,x,y}(s)=\int_{t}^{s}f^{\top }(s^{\prime },X_{t,x}(s^{\prime
}))Y_{t,x,y}(s^{\prime })ds^{\prime }+u^{\top }(s,X_{t,x}(s))Y_{t,x,y}(s).
\label{PR5}
\end{equation}%
Using Ito's formula, we get%
\begin{equation}
d\xi =\sum_{k=1}^{m}f^{k}Y^{k}ds+d(u^{\top }Y),  \label{PR05}
\end{equation}%
\begin{eqnarray}
d(u^{\top }Y) &=&\sum_{k=1}^{m}d[u^{k}Y^{k}]=\sum_{k=1}^{m}\frac{\partial
u^{k}}{\partial s}Y^{k}ds+\sum_{k=1}^{m}\sum_{i=1}^{n}\frac{\partial u^{k}}{%
\partial x^{i}}[a^{i}ds+\sum_{r=1}^{l}\sigma _{r}^{i}dw_{r}(s)]Y^{k}
\label{PR6} \\
&&+\frac{1}{2}\sum_{k=1}^{m}\sum_{i,j=1}^{n}\frac{\partial ^{2}u^{k}}{%
\partial x^{i}\partial x^{j}}\sum_{r=1}^{l}\sigma _{r}^{i}\sigma
_{r}^{j}ds\cdot Y^{k}+\sum_{k=1}^{m}u^{k}[BY]^{k}ds  \notag \\
&&+\sum_{k=1}^{m}u^{k}\left[ \sum_{r=1}^{l}\vartheta _{r}Ydw_{r}(s)\right]
^{k}+\sum_{k=1}^{m}\sum_{i=1}^{n}\frac{\partial u^{k}}{\partial x^{i}}%
dX^{i}dY^{k}.  \notag
\end{eqnarray}%
Further,%
\begin{equation}
\sum_{k=1}^{m}u^{k}[BY]^{k}ds=\sum_{k=1}^{m}[B^{\top }u]^{k}Y^{k},
\label{PR7}
\end{equation}%
\begin{eqnarray}
\sum_{k=1}^{m}\sum_{i=1}^{n}\frac{\partial u^{k}}{\partial x^{i}}%
dX^{i}dY^{k}=\sum_{k=1}^{m}\sum_{i=1}^{n}\frac{\partial u^{k}}{\partial x^{i}%
}\sum_{r=1}^{l}\sigma _{r}^{i}dw_{r}(s)\sum_{r^{\prime
}=1}^{l}\sum_{j=1}^{m}\vartheta _{r^{\prime }}^{kj}Y^{j}dw_{r^{\prime }}(s)
\label{PR8} \\
=\sum_{k=1}^{m}\sum_{i=1}^{n}\frac{\partial u^{k}}{\partial x^{i}}%
\sum_{r=1}^{l}\sum_{j=1}^{m}\sigma _{r}^{i}\vartheta
_{r}^{kj}Y^{j}ds=\sum_{k=1}^{m}(\sum_{i=1}^{n}\sum_{j=1}^{m}\sum_{r=1}^{l}%
\sigma _{r}^{i}\vartheta _{r}^{jk}\frac{\partial u^{j}}{\partial x^{i}}%
)Y^{k}ds.  \notag
\end{eqnarray}%
In (\ref{PR05})-(\ref{PR8}) all the coefficients and functions have $s,$ $%
X_{t,x}(s)$ as their arguments.

Substituting (\ref{PR6})-(\ref{PR8}) in (\ref{PR05}) and taking into account
that $u$ is a solution of (\ref{PR1}), we get%
\begin{equation}
d\xi =\sum_{r=1}^{l}\sum_{k=1}^{m}(u^{k}(\vartheta _{r}Y)^{k}+\sum_{i=1}^{n}%
\frac{\partial u^{k}}{\partial x^{i}}\sigma _{r}^{i}Y^{k})dw_{r}(s).
\label{PR9}
\end{equation}%
It is known that if
\begin{equation}
E\int_{t}^{T}\sum_{r=1}^{l}\left[ \sum_{k=1}^{m}(u^{k}(\vartheta
_{r}Y)^{k}+\sum_{i=1}^{n}\frac{\partial u^{k}}{\partial x^{i}}\sigma
_{r}^{i}Y^{k})\right] ^{2}ds<\infty  \label{PR10}
\end{equation}%
then%
\begin{equation}
E\int_{t}^{T}\sum_{r=1}^{l}\sum_{k=1}^{m}(u^{k}(\vartheta
_{r}Y)^{k}+\sum_{i=1}^{n}\frac{\partial u^{k}}{\partial x^{i}}\sigma
_{r}^{i}Y^{k})dw_{r}(s)=E[\xi _{t,x,y}(T)-\xi _{t,x,y}(t)]=0.  \label{PR11}
\end{equation}%
At the same time (see (\ref{PR5}))%
\begin{equation}
\xi _{t,x,y}(T)-\xi _{t,x,y}(t)=\int_{t}^{T}f^{\top }(s^{\prime
},X_{t,x}(s^{\prime }))Y_{t,x,y}(s^{\prime })ds^{\prime }+\varphi ^{\top
}(X_{t,x}(T))Y_{t,x,y}(T)-u^{\top }(t,x)y.  \label{PR12}
\end{equation}%
Hence%
\begin{equation}
u^{\top }(t,x)y=E\left[ \int_{t}^{T}f^{\top }(s^{\prime },X_{t,x}(s^{\prime
}))Y_{t,x,y}(s^{\prime })ds^{\prime }+\varphi ^{\top
}(X_{t,x}(T))Y_{t,x,y}(T)\right] .  \label{PR13}
\end{equation}

So, we have obtained that under certain conditions ensuring existence of a
sufficiently smooth solution of (\ref{PR1})-(\ref{PR2}), existence and
uniqueness of solution of (\ref{PR3})-(\ref{PR4}), and boundedness (\ref%
{PR10}), the probabilistic representation of the solution to the problem (%
\ref{PR1})-(\ref{PR2}) is given by formula (\ref{PR13}).

\subsection{A family of probabilistic representations}

Now we restrict ourselves to the case $\vartheta _{r}=0,\ r=1,\ldots ,l$,
i.e. (see (\ref{PR1})-(\ref{PR4})):%
\begin{equation}
\frac{\partial u^{k}}{\partial t}+\frac{1}{2}\sum_{i,j=1}^{n}\left[
\sum_{r=1}^{l}\sigma _{r}\sigma _{r}^{\top }\right] ^{ij}\frac{\partial
^{2}u^{k}}{\partial x^{i}\partial x^{j}}+\sum_{i=1}^{n}a^{i}\frac{\partial
u^{k}}{\partial x^{i}}+[B^{\top }u]^{k}+f^{k}=0,\ k=1,\ldots ,m,
\label{PR14}
\end{equation}%
\begin{equation}
u(T,x)=\varphi (x),  \label{PR15}
\end{equation}%
and%
\begin{eqnarray}
dX &=&a(s,X)ds+\sum_{r=1}^{l}\sigma _{r}(s,X)dw_{r}(s),\ X(t)=x,
\label{PR16} \\
dY &=&B(s,X)Yds,\ Y(t)=y.  \label{PR17}
\end{eqnarray}

Introduce the system
\begin{equation}
dX=a(s,X)ds-\sum_{r=1}^{l}\mu _{r}(s,X)\sigma
_{r}(s,X)ds+\sum_{r=1}^{l}\sigma _{r}(s,X)dw_{r}(s),\ X(t)=x,  \label{PR18}
\end{equation}%
\begin{equation}
dY=B(s,X)Yds,\ Y(t)=y,  \label{PR19}
\end{equation}%
\begin{equation}
dQ=\sum_{r=1}^{l}\mu _{r}(s,X)Qdw_{r}(s),\ Q(t)=1,  \label{PR20}
\end{equation}%
\begin{equation}
dZ=Qf^{\top }(s,X)Yds+Q\sum_{r=1}^{l}F_{r}^{\top }(s,X)Ydw_{r}(s),\ Z(t)=0.
\label{PR21}
\end{equation}%
In (\ref{PR18})-(\ref{PR21}) $\mu _{r},\ Q,$ and $Z$ are scalars; $\ F_{r},\
r=1,\ldots ,l,$ are column-vectors of dimension $m;\ \mu _{r}$ and $F_{r}$
are arbitrary functions, however, with good analytical properties.

Introduce the process%
\begin{equation}
\eta _{t,x,y}(s)=Q_{t,x,y,1}(s)u^{\top
}(s,X_{t,x}(s))Y_{t,x,y}(s)+Z_{t,x,y,1,0}(s).  \label{PR22}
\end{equation}%
Using Ito's formula and taking into account that $u$ is a solution of (\ref%
{PR14}), we get that under arbitrary $\mu _{r}$ and $F_{r}$%
\begin{equation}
d\eta =\sum_{r=1}^{l}Q\mathbb{\cdot (}\sum_{i=1}^{n}\frac{\partial u^{\top }%
}{\partial x^{i}}\sigma _{r}^{i}+\mu _{r}u^{\top }+F_{r}^{\top })Ydw_{r}(s).
\label{PR23}
\end{equation}

If
\begin{equation}
E\int_{t}^{T}\sum_{r=1}^{l}\left[ Q\mathbb{\cdot (}\sum_{i=1}^{n}\frac{%
\partial u^{\top }}{\partial x^{i}}\sigma _{r}^{i}+\mu _{r}u^{\top
}+F_{r}^{\top })Y\right] ^{2}ds<\infty   \label{PR24}
\end{equation}%
then%
\begin{equation}
E\int_{t}^{T}\sum_{r=1}^{l}Q\mathbb{\cdot (}\sum_{i=1}^{n}\frac{\partial
u^{\top }}{\partial x^{i}}\sigma _{r}^{i}+\mu _{r}u^{\top }+F_{r}^{\top
})Ydw_{r}(s)=E(\eta _{t,x,y}(T)-\eta _{t,x,y}(t))=0.  \label{PR25}
\end{equation}%
We have%
\begin{eqnarray}
\eta _{t,x,y}(T)-\eta _{t,x,y}(t) &=&Q(T)\varphi ^{\top
}(X_{t,x}(T))Y(T)-u(t,x)y  \label{PR26} \\
&&+\int_{t}^{T}Q_{t,x,y,1}(s^{\prime })f^{\top }(s^{\prime
},X_{t,x}(s^{\prime }))Y_{t,x,y}(s^{\prime })ds^{\prime }  \notag \\
&&+\int_{t}^{T}Q_{t,x,y,1}(s^{\prime })\sum_{r=1}^{l}F_{r}^{\top }(s^{\prime
},X_{t,x}(s^{\prime }))Y_{t,x,y}(s^{\prime })dw_{r}(s^{\prime }).  \notag
\end{eqnarray}%
Under the natural assumption%
\begin{equation*}
E\left[ \int_{t}^{T}Q_{t,x,y,1}(s^{\prime })\sum_{r=1}^{l}F_{r}^{\top
}(s^{\prime },X_{t,x}(s^{\prime }))Y_{t,x,y}(s^{\prime })dw_{r}(s^{\prime })%
\right] =0,
\end{equation*}%
using (\ref{PR25}) and (\ref{PR26}), we obtain the family of probabilistic
representations for the solution of (\ref{PR14})-(\ref{PR15}):
\begin{equation}
u^{\top }(t,x)y=E\left[ Q_{t,x,y,1}(T)\varphi ^{\top
}(X_{t,x}(T))Y_{t,x,y}(T)+\int_{t}^{T}Q_{t,x,y,1}(s^{\prime })f^{\top
}(s^{\prime },X_{t,x}(s^{\prime }))Y_{t,x,y}(s^{\prime })ds^{\prime }\right]
,  \label{PR27}
\end{equation}%
where the expressions under sign $E$ depend on a choice of $\mu _{r}$ and $%
F_{r}.$ We see that the expectation of $\eta _{t,x,y}(T)$ in the right hand
side of (\ref{PR27}) is equal to $u(t,x)y$\ and it is independent of a
choice of $\mu _{r}$ and $F_{r}.$ At the same time, the variance $Var[\eta
_{t,x,y}(T)]$ does depend on $\mu _{r}$ and $F_{r}.$

\subsection{Probabilistic representations for the vorticity}

System (\ref{V1}) has the form of (\ref{PR14}) with $m=3,\ n=3,$%
\begin{eqnarray*}
\sigma _{1} &=&\left[
\begin{array}{c}
\sigma \\
0 \\
0%
\end{array}%
\right] ,\ \sigma _{2}=\left[
\begin{array}{c}
0 \\
\sigma \\
0%
\end{array}%
\right] ,\ \sigma _{3}=\left[
\begin{array}{c}
0 \\
0 \\
\sigma%
\end{array}%
\right] , \\
a(t,x) &=&-\left[
\begin{array}{c}
u^{1}(t,x) \\
u^{2}(t,x) \\
u^{3}(t,x)%
\end{array}%
\right] ,\ B^{\top }=\left[
\begin{array}{ccc}
\partial u^{1}/\partial x^{1} & \partial u^{1}/\partial x^{2} & \partial
u^{1}/\partial x^{3} \\
\partial u^{2}/\partial x^{1} & \partial u^{2}/\partial x^{2} & \partial
u^{2}/\partial x^{3} \\
\partial u^{3}/\partial x^{1} & \partial u^{3}/\partial x^{2} & \partial
u^{3}/\partial x^{3}%
\end{array}%
\right] ,
\end{eqnarray*}%
with $\omega $ instead of $u$ and $g$ instead of $f,\ \sigma $ is a positive
constant, i.e., the more detailed writing of (\ref{V1}) has the form:
\begin{eqnarray}
\frac{\partial \omega ^{k}}{\partial t}+\frac{1}{2}\sigma ^{2}\Delta \omega
^{k}-\sum_{i=1}^{3}u^{i}(t,x)\frac{\partial \omega ^{k}}{\partial x^{i}}%
+\sum_{i=1}^{3}\frac{\partial u^{k}}{\partial x^{i}}(t,x)\omega
^{i}+g^{k}(t,x)=0,  \label{V13} \\
\omega ^{k}(T,x)=\phi ^{k}(x),\ k=1,2,3.  \label{V14}
\end{eqnarray}

Let us put $\mu _{r}(s,x)=-a^{r}(s,x),\ F^{r}(s,x)=0$ in the family of
representations (\ref{PR18})-(\ref{PR21}), (\ref{PR27}) for the problem (\ref%
{V13})-(\ref{V14}). We get%
\begin{eqnarray}
dX^{i} &=&\sigma dw_{i}(s),\ X^{i}(t)=x^{i},\ i=1,2,3,  \label{V15} \\
dY^{i} &=&\sum_{j=1}^{3}\frac{\partial u^{j}}{\partial x^{i}}(s,X)Y^{j}ds,\
Y^{i}(t)=y^{i},\ i=1,2,3,  \label{V16} \\
dQ &=&-Q\sum_{j=1}^{3}u^{j}(s,X)dw_{j}(s),\ Q(t)=1,  \label{V17} \\
dZ &=&Q\sum_{j=1}^{3}g^{j}(s,X)Y^{j}ds,\ Z(t)=0,  \label{V18} \\
\omega ^{\top }(t,x)y &=&E\left[ Q_{t,x,y,1}(T)\phi ^{\top
}(X_{t,x}(T))Y_{t,x,y}(T)+Z_{t,x,y,1,0}(T)\right] .  \label{V19}
\end{eqnarray}%
The components $\omega ^{1},\ \omega ^{2},~\omega ^{3}$ of $\omega $ are
obtained from (\ref{V19}) under $y$ equal subsequently to $(1,0,0)^{\top },\
(0,1,0)^{\top },\ (0,0,1)^{\top }.$

\begin{example}
\label{ex:31}(\textit{The Monte Carlo calculation of the Fourier coefficients%
}) Due to (\ref{vor2}), we have
\begin{eqnarray}
\omega _{\mathbf{n}}^{j}(t) &=&\frac{1}{L^{3}}\left( \omega
^{j}(t,x),e^{i(2\pi /L)(\mathbf{n},x)}\right)  \label{ex31_1} \\
&=&\frac{1}{L^{3}}\int_{Q}\omega ^{j}(t,x)e^{-i(2\pi /L)(\mathbf{n},x)}dx,\
\ \mathbf{n\in Z}^{3},\ \mathbf{n}\neq 0,\ j=1,2,3.  \notag
\end{eqnarray}%
Let $\xi $ be a random variable uniformly distributed on $Q.$ Then (\ref%
{ex31_1}) can be written as
\begin{equation*}
\omega _{\mathbf{n}}^{j}(t)=E\left[ \omega ^{j}(t,\xi )e^{-i(2\pi /L)(%
\mathbf{n},\xi )}\right] ,
\end{equation*}%
where the expectation can be approximated using the Monte Carlo technique
and hence
\begin{equation*}
\omega _{\mathbf{n}}^{j}(t)\doteq \frac{1}{M}\sum_{m=1}^{M}\omega ^{j}(t,\xi
^{(m)})e^{-i(2\pi /L)(\mathbf{n},\xi ^{(m)})}
\end{equation*}%
with $\xi ^{(m)}$ being independent realizations of $\xi .$ In turn, every $%
\omega (t,\xi ^{(m)})$ can be computed by applying the Monte Carlo technique
and weak-sense approximation of SDEs to the representation (\ref{PR27}), (%
\ref{V15})-(\ref{V19}).
\end{example}

\section{Approximation method based on vorticity\label{sec:approx}}

Let us introduce a uniform partition of the time interval $[0,T]$: $%
0=t_{0}<t_{1}<\cdots <t_{N}=T$ and the time step $h=T/N$ (we restrict
ourselves to the uniform time discretization for simplicity only). In this
section we derive an approximation for the vorticity (Section~\ref%
{sec:approxcon}) and study its properties (divergence free property in
Section~\ref{sec:approxdiv}, one-step error in Section~\ref{sec:approxerr},
and global convergence in Section~\ref{sec:conv}).

\subsection{Construction of the method\label{sec:approxcon}}

Let $\omega (t_{k+1},x),\ k=0,\ldots ,N-1,$ be known exactly. Then $%
u(t_{k+1},x)$ can be calculated exactly due to (\ref{V9}): $%
u(t_{k+1},x)=U\omega (t_{k+1},x).$ The formula (\ref{V19}),%
\begin{equation}
\omega ^{\top }(t_{k},x)y=E\left[ Q_{t_{k},x,y,1}(t_{k+1})\omega ^{\top
}(t_{k+1},X_{t_{k},x}(t_{k+1}))Y_{t_{k},x,y}(t_{k+1})+Z_{t_{k},x,y,1,0}(t_{k+1})%
\right] ,  \label{M1}
\end{equation}%
gives the value of the solution of (\ref{V1}) at $t_{k}$ assuming that $%
u(t,x),\ t_{k}\leq t<t_{k+1},$ is known exactly. We note that knowing this $%
u(t,x)$ is necessary for equations (\ref{V16})-(\ref{V17}).

Let us replace the unknown $u(t,x)$ in (\ref{V13}) by the function%
\begin{equation}
\hat{u}(t,x):=u(t_{k+1},x):=\hat{u}(x),\ t_{k}\leq t<t_{k+1}.
\label{eq:freeze}
\end{equation}%
As an approximation of $\omega (t,x)$ on $\ t_{k}\leq t\leq t_{k+1},$ we
take $\tilde{\omega}(t,x)$ satisfying the system
\begin{eqnarray}
\frac{\partial \tilde{\omega}^{i}}{\partial t}+\frac{\sigma ^{2}}{2}\Delta
\tilde{\omega}^{i}-\sum_{j=1}^{3}\hat{u}^{j}(x)\frac{\partial \tilde{\omega}%
^{i}}{\partial x^{j}}+\sum_{i=1}^{3}\frac{\partial \hat{u}^{i}}{\partial
x^{j}}(x)\tilde{\omega}^{j}+g(t,x) &=&0,\ t_{k}\leq t<t_{k+1},  \label{M2} \\
\tilde{\omega}^{i}(t_{k+1},x) &=&\omega ^{i}(t_{k+1},x),\ \ \tilde{\omega}%
^{i}(t_{k+1},x+Le_{j})=\tilde{\omega}^{i}(t_{k+1},x),\ j=1,2,3,  \label{M3}
\\
i &=&1,2,3.  \notag
\end{eqnarray}%
We observe that (\ref{M2})-(\ref{M3}) can be also obtained from (\ref{V13})
by freezing the velocity $u(t,x)$ on every time step according to (\ref%
{eq:freeze}).

Now we propose the method for solving the problem (\ref{ns1p})-(\ref{ns4p})
with negative direction of time. On the first step of the method we set
\begin{equation*}
\tilde{\omega}(t_{N},x)=\mathop{\rm curl}u(t_{N},x)=\phi (x)=\mathop{\rm
curl}\varphi (x)
\end{equation*}%
and
\begin{equation*}
\hat{u}(x)=\hat{u}(t,x)=u(t_{N},x)=\varphi (x),\ \ t_{N-1}\leq t\leq t_{N}.
\end{equation*}%
Then we solve the system (\ref{M2})-(\ref{M3}) on\textit{\ }$[t_{N-1},t_{N}]$
to obtain $\tilde{\omega}(t,x)$ and to construct
\begin{equation*}
\hat{u}(t_{N-1},x)=U\tilde{\omega}(t_{N-1},x).
\end{equation*}%
On the second step we solve (\ref{M2})-(\ref{M3}) on\textit{\ }$%
[t_{N-2},t_{N-1})$ having $\tilde{\omega}(t_{N-1},x)$ and setting $\hat{u}%
(t,x)=\hat{u}(x)=\hat{u}(t_{N-1},x)$ for $t_{N-2}\leq t<t_{N-1}.$ As a
result, we obtain $\tilde{\omega}(t,x)$ on $[t_{N-2},t_{N-1})$ and $\hat{u}%
(t_{N-2},x)=U\tilde{\omega}(t_{N-2},x),$ and so on. Proceeding in this way,
we obtain on the $N$-th step the approximation $\tilde{\omega}(t,x)$ on $%
[t_{0},t_{1})$ for $\omega (t,x)$ having $\tilde{\omega}(t_{1},x)$ and $\hat{%
u}(x)=\hat{u}(t_{1},x)=U\tilde{\omega}(t_{1},x)$ and setting $\hat{u}(t,x)=%
\hat{u}(x)=\hat{u}(t_{1},x)$ for $t_{0}\leq t<t_{1}.$ Finally, $\hat{u}%
(t_{0},x)=U\tilde{\omega}(t_{0},x).$

It is also useful to introduce
\begin{equation}
\tilde{u}(t,x):=U\tilde{\omega}(t,x),\ \ t_{0}\leq t\leq t_{N}.
\label{u_tilde}
\end{equation}%
In contrast to $\hat{u},$ the function $\tilde{u}$ is continuous in $t.$
These functions coincide at $t=t_{k},$ $k=0,\ldots ,N.$

At each step of this method one has to solve the system (\ref{M2})-(\ref{M3}%
). In contrast to the system (\ref{ns1p})-(\ref{ns4p}), the system (\ref{M2}%
)-(\ref{M3}) does not have the divergence-free condition and it is linear.
Then the solution of (\ref{M2})-(\ref{M3}) can be found using probabilistic
representations. We pay attention to the fact that in the vorticity
formulation of the NSE the pressure term disappears.

In order to realise the approximation process described above, it is
sufficient that on every time interval $[t_{k},t_{k+1}],$ $k=N-1,N-2,\ldots
,1,0,$ there exists a solution of the linear parabolic system (\ref{M2})-(%
\ref{M3}) (we denote such a solution $\tilde{\omega}_{k}(t,x))$ which
satisfies the condition
\begin{equation}
\tilde{\omega}_{k}(t_{k+1},x)=\left\{
\begin{array}{c}
\mathop{\rm curl}\varphi (x),\ k=N-1, \\
\tilde{\omega}_{k+1}(t_{k+1},x),\ k=N-2,\ldots ,0,%
\end{array}%
\right.  \label{vor45}
\end{equation}%
and has the time-independent $\hat{u}(x)$ within each interval $%
[t_{k},t_{k+1})$ $\ $defined as\
\begin{equation}
\hat{u}(x):=\hat{u}_{k}(x)=U\tilde{\omega}_{k}(t_{k+1},x),\ t_{k}\leq
t<t_{k+1}.  \label{vor46}
\end{equation}%
Clearly, $\hat{u}(x)$ used in (\ref{M2}) are different on the time intervals
$[t_{k},t_{k+1})$.

\subsection{The divergence-free property of the method\label{sec:approxdiv}}

The evolution equation (\ref{V1}) for vorticity has the form%
\begin{equation*}
\frac{\partial \omega }{\partial t}=\mathop{\rm curl}[\ldots ].
\end{equation*}%
Due to this fact, any solution of (\ref{V1}) with $\mathop{\rm div}\omega
(t_{k},x)=0$ is divergence free for $t\leq t_{k}:$ $\mathop{\rm div}\omega
(t,x)=0, $ $t\leq t_{k}.$ Indeed, this property can be seen after applying
the operator $\mathop{\rm div}$ to (\ref{V1}) and taking into account the
equality $\mathop{\rm div}\mathop{\rm curl}[\ldots ]=0.$

A very important property of the proposed method is that the constructed
approximation $\tilde{\omega}_{k}(t,x)$ is also divergent free.

\begin{theorem}
\label{thm:divfree}The solution $\tilde{\omega}(t,x),\ t_{k}\leq t\leq
t_{k+1},$\ of $(\ref{M2})$-$(\ref{M3})$ is divergent free.
\end{theorem}

\noindent \textbf{Proof}. Let us take $\mathop{\rm div}$\ of the equation (%
\ref{M2}). In (\ref{M2}) we have that $\hat{u}(t,x)=\hat{u}%
(x)=u(t_{k+1},x),\ t_{k}\leq t<t_{k+1},$ and $\hat{u}(x)$\ is\ divergent
free: $\mathop{\rm div}\hat{u}=0$. Besides, $\mathop{\rm div}g=0$. We have%
\begin{equation*}
-(\hat{u},\nabla )\tilde{\omega}=-(\hat{u}^{1}\frac{\partial }{\partial x^{1}%
}+\hat{u}^{2}\frac{\partial }{\partial x^{2}}+\hat{u}^{3}\frac{\partial }{%
\partial x^{3}})\tilde{\omega}=-\left[
\begin{array}{c}
\hat{u}^{1}\frac{\partial \tilde{\omega}^{1}}{\partial x^{1}}+\hat{u}^{2}%
\frac{\partial \tilde{\omega}^{1}}{\partial x^{2}}+\hat{u}^{3}\frac{\partial
\tilde{\omega}^{1}}{\partial x^{3}} \\
\hat{u}^{1}\frac{\partial \tilde{\omega}^{2}}{\partial x^{1}}+\hat{u}^{2}%
\frac{\partial \tilde{\omega}^{2}}{\partial x^{2}}+\hat{u}^{3}\frac{\partial
\tilde{\omega}^{2}}{\partial x^{3}} \\
\hat{u}^{1}\frac{\partial \tilde{\omega}^{3}}{\partial x^{1}}+\hat{u}^{2}%
\frac{\partial \tilde{\omega}^{3}}{\partial x^{2}}+\hat{u}^{3}\frac{\partial
\tilde{\omega}^{3}}{\partial x^{3}}%
\end{array}%
\right] ,
\end{equation*}%
\begin{eqnarray*}
\mathop{\rm div}[-(\hat{u},\nabla )\tilde{\omega}]=-(\frac{\partial \hat{u}%
^{1}}{\partial x^{1}}\frac{\partial \tilde{\omega}^{1}}{\partial x^{1}}+%
\frac{\partial \hat{u}^{2}}{\partial x^{1}}\frac{\partial \tilde{\omega}^{1}%
}{\partial x^{2}}+\frac{\partial \hat{u}^{3}}{\partial x^{1}}\frac{\partial
\tilde{\omega}^{1}}{\partial x^{3}} \\
\frac{\partial \hat{u}^{1}}{\partial x^{2}}\frac{\partial \tilde{\omega}^{2}%
}{\partial x^{1}}+\frac{\partial \hat{u}^{2}}{\partial x^{2}}\frac{\partial
\tilde{\omega}^{2}}{\partial x^{2}}+\frac{\partial \hat{u}^{3}}{\partial
x^{2}}\frac{\partial \tilde{\omega}^{2}}{\partial x^{3}}+\frac{\partial \hat{%
u}^{1}}{\partial x^{3}}\frac{\partial \tilde{\omega}^{3}}{\partial x^{1}}+%
\frac{\partial \hat{u}^{2}}{\partial x^{3}}\frac{\partial \tilde{\omega}^{3}%
}{\partial x^{2}}+\frac{\partial \hat{u}^{3}}{\partial x^{3}}\frac{\partial
\tilde{\omega}^{3}}{\partial x^{3}}) \\
-(\hat{u}^{1}\frac{\partial }{\partial x^{1}}\mathop{\rm div}\tilde{\omega}+%
\hat{u}^{2}\frac{\partial }{\partial x^{2}}\mathop{\rm div}\tilde{\omega}+%
\hat{u}^{3}\frac{\partial }{\partial x^{3}}\mathop{\rm div}\tilde{\omega}).
\end{eqnarray*}%
Analogously,%
\begin{eqnarray*}
\mathop{\rm div}[(\tilde{\omega},\nabla )\hat{u}]=\frac{\partial \tilde{%
\omega}^{1}}{\partial x^{1}}\frac{\partial \hat{u}^{1}}{\partial x^{1}}+%
\frac{\partial \tilde{\omega}^{2}}{\partial x^{1}}\frac{\partial \hat{u}^{1}%
}{\partial x^{2}}+\frac{\partial \tilde{\omega}^{3}}{\partial x^{1}}\frac{%
\partial \hat{u}^{1}}{\partial x^{3}} \\
\frac{\partial \tilde{\omega}^{1}}{\partial x^{2}}\frac{\partial \hat{u}^{2}%
}{\partial x^{1}}+\frac{\partial \tilde{\omega}^{2}}{\partial x^{2}}\frac{%
\partial \hat{u}^{2}}{\partial x^{2}}+\frac{\partial \tilde{\omega}^{3}}{%
\partial x^{2}}\frac{\partial \hat{u}^{2}}{\partial x^{3}}+\frac{\partial
\tilde{\omega}^{1}}{\partial x^{3}}\frac{\partial \hat{u}^{3}}{\partial x^{1}%
}+\frac{\partial \tilde{\omega}^{2}}{\partial x^{3}}\frac{\partial \hat{u}%
^{3}}{\partial x^{2}}+\frac{\partial \tilde{\omega}^{3}}{\partial x^{3}}%
\frac{\partial \hat{u}^{3}}{\partial x^{3}} \\
+(\tilde{\omega}^{1}\frac{\partial }{\partial x^{1}}\mathop{\rm div}\hat{u}+%
\tilde{\omega}^{2}\frac{\partial }{\partial x^{2}}\mathop{\rm div}\hat{u}+%
\tilde{\omega}^{3}\frac{\partial }{\partial x^{3}}\mathop{\rm div}\hat{u}).
\end{eqnarray*}%
Since $\mathop{\rm div}\hat{u}=0,$ we get%
\begin{eqnarray*}
&&\mathop{\rm div}[-(\hat{u},\nabla )\tilde{\omega}]+\mathop{\rm
div}[(\tilde{\omega},\nabla )\hat{u}] \\
&=&-(\hat{u}^{1}\frac{\partial }{\partial x^{1}}\mathop{\rm div}\tilde{\omega%
}+\hat{u}^{2}\frac{\partial }{\partial x^{2}}\mathop{\rm div}\tilde{\omega}+%
\hat{u}^{3}\frac{\partial }{\partial x^{3}}\mathop{\rm div}\tilde{\omega}).
\end{eqnarray*}%
Hence, taking $\mathop{\rm div}$ of (\ref{M2}) gives the following equation
for $\mathop{\rm div}\tilde{\omega}:$%
\begin{eqnarray}
\frac{\partial \mathop{\rm div}\tilde{\omega}}{\partial t}-(\hat{u}^{1}\frac{%
\partial }{\partial x^{1}}\mathop{\rm div}\tilde{\omega}+\hat{u}^{2}\frac{%
\partial }{\partial x^{2}}\mathop{\rm div}\tilde{\omega}+\hat{u}^{3}\frac{%
\partial }{\partial x^{3}}\mathop{\rm div}\tilde{\omega})+\frac{\sigma ^{2}}{%
2}\Delta \mathop{\rm div}\tilde{\omega}=0,  \label{M4} \\
t_{k}\leq t<t_{k+1},\ \ \mathop{\rm div}\tilde{\omega}(t_{k+1},x)=0.
\label{M5}
\end{eqnarray}%
From here, due to uniqueness of solution to the problem (\ref{M4})-(\ref{M5}%
), we obtain%
\begin{equation*}
\mathop{\rm div}\tilde{\omega}(t,x)=0,\ t_{k}\leq t\leq t_{k+1},\ x\in
\mathbf{R}^{3}.
\end{equation*}%
Theorem~\ref{thm:divfree} is proved.

\subsection{The one-step error of the method\label{sec:approxerr}}

For estimating the local error (the one-step error) in the 2D case, together
with the solution $\omega (t,x),\ t_{k}\leq t\leq t_{k+1},$ of (\ref{M6}),
we consider the approximation $\tilde{\omega}(t,x),\ $which satisfies the
equation%
\begin{equation}
\frac{\partial \tilde{\omega}}{\partial t}-\hat{u}^{1}(x)\frac{\partial
\tilde{\omega}}{\partial x^{1}}(t,x)-\hat{u}^{2}(x)\frac{\partial \tilde{%
\omega}}{\partial x^{2}}(t,x)+\frac{\sigma ^{2}}{2}\Delta \tilde{\omega}%
(t,x)+g(t,x^{1},x^{2})=0  \label{M7}
\end{equation}%
and the Cauchy condition%
\begin{equation}
\tilde{\omega}(t_{k+1},x)=\omega (t_{k+1},x).  \label{M07}
\end{equation}%
The difference
\begin{equation*}
\delta _{\omega }(t,x):=\omega (t,x)-\tilde{\omega}(t,x),
\end{equation*}%
which is the one step error, is a solution to the problem%
\begin{gather}
\frac{\partial \delta _{\omega }}{\partial t}+\frac{\sigma ^{2}}{2}\Delta
\delta _{\omega }-u^{1}\frac{\partial \delta _{\omega }}{\partial x^{1}}%
-u^{2}\frac{\partial \delta _{\omega }}{\partial x^{2}}-(u^{1}-\hat{u}^{1})%
\frac{\partial \tilde{\omega}}{\partial x^{1}}-(u^{2}-\hat{u}^{2})\frac{%
\partial \tilde{\omega}}{\partial x^{2}}=0,  \label{M8} \\
\delta _{\omega }(t_{k+1},x)=0.  \label{M9}
\end{gather}

\begin{theorem}
\label{thm:onestep2D} The one-step error of $\tilde{\omega}(t,x),\ t_{k}\leq
t\leq t_{k+1},$\ which solves $(\ref{M2})$-$(\ref{M3})$ is of second order
with respect to $h:$%
\begin{equation}
|\delta _{\omega }(t,x)|\leq Kh^{2},\ t_{k}\leq t\leq t_{k+1},\ x\in \mathbf{%
R}^{2}.  \label{M10}
\end{equation}
\end{theorem}

\noindent \textbf{Proof}. Let us write the probabilistic representation of
the form (\ref{V15})-(\ref{V19}) for the solution to problem (\ref{M8})-(\ref%
{M9}):%
\begin{equation}
dX^{i}=\sigma dw_{i}(s),\ X^{i}(t)=x^{i},\ i=1,2,  \label{M11}
\end{equation}%
\begin{equation}
dQ=-Q(u^{1}dw_{1}+u^{2}dw_{2}),\ Q(t)=1,  \label{M12}
\end{equation}%
\begin{equation}
dZ=-Q((u^{1}-\hat{u}^{1})\frac{\partial \tilde{\omega}}{\partial x^{1}}%
+(u^{2}-\hat{u}^{2})\frac{\partial \tilde{\omega}}{\partial x^{2}})ds,\
Z(t)=0,  \label{M13}
\end{equation}%
\begin{equation}
\delta _{\omega }(t,x)=-E\int_{t}^{t_{k+1}}Q((u^{1}-\hat{u}^{1})\frac{%
\partial \tilde{\omega}}{\partial x^{1}}+(u^{2}-\hat{u}^{2})\frac{\partial
\tilde{\omega}}{\partial x^{2}})ds.  \label{M14}
\end{equation}%
Using boundedness of $\partial \tilde{\omega}/\partial x^{i},\ i=1,2,$ and
the inequalities%
\begin{equation*}
|u^{i}(s,X_{t,x}(s))-\hat{u}%
^{i}(X_{t,x}(s))|=|u^{i}(s,X_{t,x}(s))-u^{i}(t_{k+1},X_{t,x}(s))|\leq Ch,
\end{equation*}%
for $t_{k}\leq s\leq t_{k+1},$ we get%
\begin{equation*}
|\delta _{\omega }(t,x)|\leq \int_{t}^{t_{k+1}}E|Q|ds\cdot Kh.
\end{equation*}%
But $Q>0$ and%
\begin{equation*}
E|Q|=EQ=1,
\end{equation*}%
whence (\ref{M10}) follows. Theorem~\ref{thm:onestep2D} is proved.\bigskip

Introduce the one-step error for $\tilde{u}(t,x)$ from (\ref{u_tilde}):
\begin{equation}
\delta _{u}(t,x):=u(t,x)-U\tilde{\omega}(t,x)=U\delta _{\omega }(t,x),
\label{delta_u}
\end{equation}%
where $\tilde{\omega}(t,x),\ t_{k}\leq t\leq t_{k+1},$\ is the solution of (%
\ref{M2})-(\ref{M3}).

\begin{corollary}
\label{Cor2D}The one-step error of $\tilde{u}(t,x)$ from (\ref{u_tilde}) is
of second order with respect to $h$ in $L^{2}$-norm:%
\begin{equation}
||\delta _{u}(t,x)||_{L^{2}}\leq Kh^{2},\ t_{k}\leq t\leq t_{k+1}.
\label{Cor2D1}
\end{equation}
\end{corollary}

\noindent \textbf{Proof}. Let the Fourier coefficients for $\delta _{\omega
}(t,x)$ be $(\delta _{\omega }(t,\cdot ))_{\mathbf{n}},$ i.e.
\begin{equation*}
\delta _{\omega }(t,x)=\sum (\delta _{\omega }(t,\cdot ))_{\mathbf{n}%
}e^{i(2\pi /L)(\mathbf{n},x)}.
\end{equation*}%
Hence (cf. (\ref{Mnew1}))%
\begin{equation*}
\delta _{u}(t,x)=\frac{iL}{2\pi }\sum \frac{1}{|\mathbf{n}|^{2}}e^{i(2\pi
/L)(\mathbf{n},x)}(\delta _{\omega }(t,\cdot ))_{\mathbf{n}}\left[
\begin{array}{c}
n^{2} \\
-n^{1}%
\end{array}%
\right] ,
\end{equation*}%
i.e., the Fourier coefficients for $\delta _{u}(t,x)$ are
\begin{equation*}
(\delta _{u}(t,\cdot ))_{\mathbf{n}}=\frac{iL}{2\pi }\frac{1}{|\mathbf{n}%
|^{2}}(\delta _{\omega }(t,\cdot ))_{\mathbf{n}}\left[
\begin{array}{c}
n^{2} \\
-n^{1}%
\end{array}%
\right] .
\end{equation*}%
Then, by Parseval's identity (\ref{eq:pars}), we have
\begin{eqnarray}
||\delta _{u}(t,\cdot )||_{L^{2}} &=&\int_{Q}|\delta
_{u}(t,x)|^{2}dx=L^{2}\sum |(\delta _{u}(t,\cdot ))_{\mathbf{n}}|^{2}
\label{extra} \\
&=&\frac{L^{4}}{4\pi ^{2}}\sum |(\delta _{\omega }(t,\cdot ))_{\mathbf{n}%
}|^{2}\frac{\left( n^{1}\right) ^{2}+\left( n^{2}\right) ^{2}}{|\mathbf{n}%
|^{4}}  \notag \\
&=&\frac{L^{4}}{4\pi ^{2}}\sum \frac{|(\delta _{\omega }(t,\cdot ))_{\mathbf{%
n}}|^{2}}{|\mathbf{n}|^{2}}\leq \frac{L^{4}}{4\pi ^{2}}\sum |(\delta
_{\omega }(t,\cdot ))_{\mathbf{n}}|^{2}  \notag \\
&=&\frac{L^{2}}{4\pi ^{2}}\int_{Q}|\delta _{\omega }(t,x)|^{2}dx\leq \frac{%
L^{2}}{4\pi ^{2}}\max_{x}|\delta _{\omega }(t,x)|^{2},  \notag
\end{eqnarray}%
which together with (\ref{M10}) implies (\ref{Cor2D1}). Corollary~\ref{Cor2D}
is proved.\bigskip

The result of Theorem~\ref{thm:onestep2D} is carried over to the 3D case
without any substantial changes in the proof. In the 3D case the difference $%
\delta _{\omega }(t,x):=\omega (t,x)-\tilde{\omega}(t,x)$ is a solution to
the problem%
\begin{gather}
\frac{\partial \delta _{\omega }}{\partial t}+\frac{\sigma ^{2}}{2}\Delta
\delta _{\omega }-\sum_{i=1}^{3}u^{i}\frac{\partial \delta _{\omega }}{%
\partial x^{i}}+\sum_{i=1}^{3}\frac{\partial u}{\partial x^{i}}\delta
_{\omega }^{i}-\sum_{i=1}^{3}(u^{i}-\hat{u}^{i})\frac{\partial \tilde{\omega}%
}{\partial x^{i}}+\sum_{i=1}^{3}(\frac{\partial u}{\partial x^{i}}-\frac{%
\partial \hat{u}}{\partial x^{i}})\tilde{\omega}^{i}=0,  \label{M15} \\
\delta _{\omega }(t_{k+1},x)=0.  \label{M16}
\end{gather}

\begin{theorem}
\label{thm:onestep3D}The one-step error of $\tilde{\omega}(t,x),\ t_{k}\leq
t\leq t_{k+1},$\ which solves $(\ref{M2})$-$(\ref{M3}),$ is of second order
with respect to $h:$%
\begin{equation}
|\delta _{\omega }(t,x)|\leq Kh^{2},\ t_{k}\leq t\leq t_{k+1},\ x\in \mathbf{%
R}^{3}.  \label{M017}
\end{equation}
\end{theorem}

\noindent \textbf{Proof}. We apply the probabilistic representation (\ref%
{V15})-(\ref{V19}) to the solution of (\ref{M15})-(\ref{M16}):%
\begin{equation}
dX^{i}=\sigma dw_{i}(s),\ X^{i}(t)=x^{i},\ i=1,2,3,  \label{M17}
\end{equation}%
\begin{equation}
dY^{i}=\sum_{j=1}^{3}\frac{\partial u^{j}}{\partial x^{i}}Y^{j}ds,\
Y^{i}(t)=y^{i},\ i=1,2,3,  \label{M18}
\end{equation}%
\begin{equation}
dQ=-Q\sum_{j=1}^{3}u^{j}(s,X)dw_{j}(s),\ Q(t)=1,  \label{M19}
\end{equation}%
\begin{equation}
dZ=Q\left[ \sum_{j=1}^{3}\sum_{i=1}^{3}(\frac{\partial u^{j}}{\partial x^{i}}%
-\frac{\partial \hat{u}^{j}}{\partial x^{i}})\tilde{\omega}%
^{i}Y^{j}-\sum_{j=1}^{3}\sum_{i=1}^{3}(u^{i}-\hat{u}^{i})\frac{\partial
\tilde{\omega}^{j}}{\partial x^{i}}Y^{j}\right] ds,\ Z(t)=0,  \label{M20}
\end{equation}%
\begin{eqnarray}
\delta _{\omega }^{\top }(t,x)y=EZ_{t,x,y,1,0}(t_{k+1})  \label{M21} \\
=E\int_{t}^{t_{k+1}}Q_{t,x,y,1}(s)\left[ \sum_{j=1}^{3}\sum_{i=1}^{3}(\frac{%
\partial u^{j}}{\partial x^{i}}(s,X_{t,x}(s))-\frac{\partial \hat{u}^{j}}{%
\partial x^{i}}(X_{t,x}(s)))\tilde{\omega}^{i}(s,X_{t,x}(s))Y^{j}(s)\right.
\notag \\
-\left. \sum_{j=1}^{3}\sum_{i=1}^{3}(u^{i}(s,X_{t,x}(s))-\hat{u}%
^{i}(X_{t,x}(s)))\frac{\partial \tilde{\omega}^{j}}{\partial x^{i}}%
(s,X_{t,x}(s))Y^{j}(s)\right] ds.  \notag
\end{eqnarray}%
Using boundedness of $\tilde{\omega}^{i},\ \partial \tilde{\omega}%
^{j}/\partial x^{i},\ Y^{i}(s),\ i,j=1,2,3,$\ the inequalities%
\begin{equation*}
|u^{i}(s,X_{t,x}(s))-\hat{u}%
^{i}(X_{t,x}(s))|=|u^{i}(s,X_{t,x}(s))-u^{i}(t_{k+1},X_{t,x}(s))|\leq Ch,
\end{equation*}%
\begin{equation*}
|\frac{\partial u^{j}}{\partial x^{i}}(s,X_{t,x}(s))-\frac{\partial \hat{u}%
^{j}}{\partial x^{i}}(X_{t,x}(s))|=|\frac{\partial u^{j}}{\partial x^{i}}%
(s,X_{t,x}(s))-\frac{\partial u^{j}}{\partial x^{i}}(t_{k+1},X_{t,x}(s))|%
\leq Ch,
\end{equation*}%
for $t_{k}\leq s<t_{k+1},$ and the properties $Q>0,\ E|Q|=EQ=1,$ we get (\ref%
{M017}). Theorem~\ref{thm:onestep3D} is proved.\bigskip

We note that the one-step error estimate (\ref{Cor2D1}) for $\tilde{u}$ from
Corollary~\ref{Cor2D} is also valid in the three-dimensional case.

\subsection{Convergence theorems\label{sec:conv}}

In this section we first consider the global error for the approximation $%
\tilde{u}(t,x)$ from (\ref{u_tilde}), i.e., we are interested in estimating
the difference
\begin{equation*}
D_{\tilde{u}}:=u(t_{0},x)-\tilde{u}(t_{0},x),
\end{equation*}%
where $u(t_{0},x)$ is the solution of the NSE (\ref{ns1p})-(\ref{ns4p}).

Let us introduce the auxiliary functions $_{k}u(t,x)$ on the time intervals $%
[t_{0},t_{k}],$ $k=1,\ldots ,N:$
\begin{equation}
_{k}u(t,x):=u(t,x;t_{k},\tilde{u}(t_{k},\cdot )),\ t_{0}\leq t\leq t_{k},
\label{eq:k_u}
\end{equation}%
where $u(t,x;t_{k},\tilde{u}(t_{k},\cdot ))$ denotes the solution of the NSE
(\ref{ns1p})-(\ref{ns4p}) with the terminal condition $\varphi (\cdot )=%
\tilde{u}(t_{k},\cdot )$ prescribed at $T=t_{k}$. To prove the convergence
theorem, we assume that all the functions $_{k}u(t,x)$ are bounded together
with their derivatives up to some order.

Since $\tilde{u}(t_{N},x)=u(t_{N},x),$ we have $_{N}u(t,x)=u(t,x),$ $%
t_{0}\leq t\leq t_{N}.$ Also, note that $\tilde{u}(t_{0},x)=\
_{0}u(t_{0},x). $ Then we can re-write the global error as%
\begin{equation}
D_{\tilde{u}}=\sum_{k=0}^{N-1}\left( \ _{k+1}u(t_{0},x)-\
_{k}u(t_{0},x)\right) .  \label{c01}
\end{equation}%
We have
\begin{eqnarray}
_{k+1}u(t_{0},x) &=&u(t_{0},x;t_{k+1},\tilde{u}(t_{k+1},\cdot
))=u(t_{0},x;t_{k},u(t_{k},\cdot ;t_{k+1},\tilde{u}(t_{k+1},\cdot ))),
\label{c02} \\
_{k}u(t_{0},x) &=&u(t_{0},x;t_{k},\tilde{u}(t_{k},\cdot )).  \notag
\end{eqnarray}%
Note that the difference
\begin{equation*}
_{k}\delta _{u}(t_{k},x)=u(t_{k},x;t_{k+1},\tilde{u}(t_{k+1},\cdot ))-\tilde{%
u}(t_{k},x)
\end{equation*}%
is a one-step error (see (\ref{delta_u})), which $L^{2}$-estimate is of
order $h^{2}$ according to Corollary~\ref{Cor2D}. We remark that $%
_{k+1}u(t_{0},x)-\ _{k}u(t_{0},x)$ is the propagation error which is due to
the error in the terminal condition propagated along the trajectory of the
NSE solution.

To estimate the propagation error, we are making use of the basic energy
estimate from \cite[p. 89]{MB}, where it is proved in the whole space, but
it can be derived for the periodic case as well. In our case this energy
estimate takes the form%
\begin{equation}
\sup_{t_{0}\leq t\leq t_{k}}||\ _{k+1}u(t_{0},\cdot )-\ _{k}u(t_{0},\cdot
)||_{L^{2}}\leq C||\ _{k+1}u(t_{k},\cdot )-\ _{k}u(t_{k},\cdot )||_{L^{2}},
\label{c03}
\end{equation}%
where the constant $C>0$ depends on the function $_{k}u(t,x).$

Due to (\ref{Cor2D1}) and (\ref{c03}), we obtain
\begin{equation}
||\ _{k+1}u(t_{0},\cdot )-\ _{k}u(t_{0},\cdot )||_{L^{2}}\leq Kh^{2},
\label{c04}
\end{equation}%
where $K>0$ combines the constant $K$ from (\ref{Cor2D1}) and $C$ from (\ref%
{c03}). From (\ref{c04}) and (\ref{c01}), we get
\begin{equation*}
||D_{\tilde{u}}||_{L^{2}}\leq Kh.
\end{equation*}%
Thus, we have proved the following theorem.

\begin{theorem}
\label{thm:conv} The approximation $\tilde{u}(t,x)$ from $(\ref{u_tilde})$
for the solution of the NSE $(\ref{ns1p})$-$(\ref{ns4p})$ is of first order
in $h.$
\end{theorem}

We note that the proof of Theorem~\ref{thm:conv} tacitly used an assumption
of existence, uniqueness and regularity of solutions of the NSE problems
involved in the error estimates. Such an assumption is natural to make in
the work aimed at deriving approximations and we do not consider here how
one can prove such properties of $_{k}u(t,x)$ from (\ref{eq:k_u}).

Now we analyse the global error of $\tilde{\omega}_{k}(t,x).$

\begin{theorem}
\label{thm:conv2} The approximation $\tilde{\omega}_{k}(t,x)$ $($see $(\ref%
{u_tilde}),$\ $(\ref{vor45}))$ for the solution of the NSE $(\ref{ns1p})$-$(%
\ref{ns4p})$ converges with order 1 in $L^{2}$-norm.
\end{theorem}

\noindent \textbf{Proof}. \ Let $D_{\tilde{\omega}}(t,x;k)$ be the global
error for $\tilde{\omega}$ on the interval $[t_{k,}t_{k+1}],$ i.e.%
\begin{equation*}
D_{\tilde{\omega}}(t,x;k):=\omega (t,x)-\tilde{\omega}_{k}(t,x),
\end{equation*}%
and $D_{\hat{u}}(t,x;k)$ be the global error for $\hat{u}_{k}$ on the
interval $[t_{k,}t_{k+1}),$ i.e.%
\begin{equation*}
D_{\hat{u}}(t,x;k):=u(t,x)-\hat{u}_{k}(x).
\end{equation*}%
We have analogously to (\ref{M8})-(\ref{M9}):%
\begin{gather}
-\frac{\partial D_{\tilde{\omega}}(t,x;k)}{\partial t}=\frac{\sigma ^{2}}{2}%
\Delta D_{\tilde{\omega}}(t,x;k)-(u(t,x),\nabla )D_{\tilde{\omega}}(t,x;k)
\notag \\
-(D_{\hat{u}}(t,x;k),\nabla )\tilde{\omega}_{k}(t,x),\ \ t_{k}\leq
t<t_{k+1},\ \ k=N-1,\ldots ,0,  \label{D1} \\
D_{\tilde{\omega}}(t_{N},x;N-1)=0,  \label{D2} \\
D_{\tilde{\omega}}(t_{k+1},x;k)=D_{\tilde{\omega}}(t_{k+1},x;k+1),\ \ \
k=N-2,\ldots ,0.  \label{D3}
\end{gather}%
Then
\begin{gather}
-\frac{1}{2}\frac{d||D_{\tilde{\omega}}(t,\cdot ;k)||^{2}}{dt}=-\frac{\sigma
^{2}}{2}||\nabla D_{\tilde{\omega}}||^{2}-((u,\nabla )D_{\tilde{\omega}},D_{%
\tilde{\omega}})-((D_{\hat{u}}(t,\cdot ;k),\nabla )\tilde{\omega}_{k},D_{%
\tilde{\omega}}),  \label{sp1} \\
t_{k}\leq t<t_{k+1}.  \notag
\end{gather}%
Since $u$ is divergence free, we get%
\begin{equation}
((u,\nabla )D_{\tilde{\omega}},D_{\tilde{\omega}})=0.  \label{sp11}
\end{equation}%
Note that (see (\ref{vor46})):
\begin{eqnarray*}
D_{\hat{u}}(t,x;k) &=&u(t,x)-\hat{u}_{k}(x)=U\omega (t,x)-U\tilde{\omega}%
_{k}(t_{k+1},x) \\
&=&UD_{\tilde{\omega}}(t,x;k)+U\tilde{\omega}_{k}(t,x)-U\tilde{\omega}%
_{k}(t_{k+1},x).
\end{eqnarray*}%
Then, using (\ref{A66}) with $c_{2}=\sigma ^{2}/2$, we obtain for some $K>0$:%
\begin{eqnarray}
|((D_{\hat{u}}(t,\cdot ;k),\nabla )\tilde{\omega}_{k},D_{\tilde{\omega}})|
&=&|((U\left( D_{\tilde{\omega}}+\tilde{\omega}_{k}(t,\cdot )-\tilde{\omega}%
_{k}(t_{k+1},\cdot )\right) ,\nabla )\tilde{\omega}_{k},D_{\tilde{\omega}})|
\label{sp2} \\
&\leq &\frac{\sigma ^{2}}{2}||\nabla D_{\tilde{\omega}}||^{2}+K||D_{\tilde{%
\omega}}||^{2}  \notag \\
&&+K||\nabla \tilde{\omega}_{k}||^{2}||D_{\tilde{\omega}}+\tilde{\omega}%
_{k}(t,\cdot )-\tilde{\omega}_{k}(t_{k+1},\cdot )||^{2}  \notag \\
&\leq &\frac{\sigma ^{2}}{2}||\nabla D_{\tilde{\omega}}||^{2}+K||D_{\tilde{%
\omega}}||^{2}+K||\nabla \tilde{\omega}_{k}||^{2}||D_{\tilde{\omega}}||^{2}
\notag \\
&&+K||\nabla \tilde{\omega}_{k}||^{2}||\tilde{\omega}_{k}(t,\cdot )-\tilde{%
\omega}_{k}(t_{k+1},\cdot )||^{2}.  \notag
\end{eqnarray}%
Using boundedness of $||\frac{d}{dt}\tilde{\omega}_{k}(t,\cdot )||^{2},$ we
get
\begin{equation*}
||\tilde{\omega}_{k}(t,\cdot )-\tilde{\omega}_{k}(t_{k+1},\cdot )||^{2}\leq
Kh^{2},
\end{equation*}%
which together with boundedness of $||\nabla \tilde{\omega}_{k}||^{2}$
implies
\begin{equation}
|((D_{\hat{u}}(t,\cdot ;k),\nabla )\tilde{\omega}_{k},D_{\tilde{\omega}%
})|\leq \frac{\sigma ^{2}}{2}||\nabla D_{\tilde{\omega}}||^{2}+K||D_{\tilde{%
\omega}}||^{2}+Kh^{2}.  \label{sp22}
\end{equation}%
It follows from (\ref{sp1}), (\ref{sp11})\ and (\ref{sp22}) that
\begin{equation*}
-\frac{d\left( ||D_{\tilde{\omega}}(t,\cdot ;k)||^{2}+h^{2}\right) }{||D_{%
\tilde{\omega}}(t,\cdot ;k)||^{2}+h^{2}}\leq 2Kdt,\ t_{k}\leq t<t_{k+1}.
\end{equation*}%
Then
\begin{equation*}
||D_{\tilde{\omega}}(t_{k},\cdot ;k)||^{2}+h^{2}\leq e^{2Kh}\left( ||D_{%
\tilde{\omega}}(t_{k+1},\cdot ,k)||^{2}+h^{2}\right) .
\end{equation*}%
From here and due to (\ref{D3}), we get
\begin{equation*}
||D_{\tilde{\omega}}(t_{k},\cdot ;k-1)||^{2}\leq e^{2Kh}||D_{\tilde{\omega}%
}(t_{k+1},\cdot ,k)||^{2}+\left( e^{2Kh}-1\right) h^{2},\ k=N-1,\ldots ,1.
\end{equation*}%
Denoting $R_{k}:=||D_{\tilde{\omega}}(t_{k+1},\cdot ;k)||^{2},$ $%
k=N-1,\ldots ,0,$ we obtain (see (\ref{D2})):
\begin{eqnarray*}
R_{k-1} &\leq &e^{2Kh}R_{k}+\left( e^{2Kh}-1\right) h^{2},\ k=N-1,\ldots ,1,
\\
R_{N-1} &=&0,
\end{eqnarray*}%
and using the discrete Gronwall lemma (see e.g. \cite[p. 7]{MT1}), we arrive
at $R_{0}=||D_{\tilde{\omega}}(t_{1},\cdot ;0)||^{2}\leq Kh^{2}.$ Theorem~%
\ref{thm:conv2} is proved.

\section{Stochastic Navier-Stokes equations\label{sec:sns}}

In this section we carry over the results of Section~\ref{sec:approx} for
the deterministic NSE to two-dimensional NSE with additive noise. After
introducing the stochastic NSE in velocity-vorticity formulation, we prove
two auxiliary lemmas (Section~\ref{sec:sns1}) about its solution; we
consider a one-step approximation of vorticity and its properties (Section~%
\ref{sec:sns2}); we introduce the numerical method for vorticity and prove
boundedness of its moments in Section~\ref{sec:sns3}; and, finally, we prove
first-order mean-square convergence of the method in Section~\ref{sec:sns4}.
The global convergence proof contains ideas, which can potentially be
exploited in analysis of numerical methods for a wider class of semilinear
SPDEs.

Let $(\Omega ,\mathcal{F},P)$ be a probability space and $(w(t),\mathcal{F}%
_{t}^{w})=((w_{1}(t),\ldots ,w_{q}(t))^{\top },\mathcal{F}_{t})$ be a $q$%
-dimensional standard Wiener process, where $\mathcal{F}_{t},\ 0\leq t\leq
T, $ is an increasing family of $\sigma $-subalgebras of $\mathcal{F}$
induced by $w(t).$ We consider the system of stochastic Navier-Stokes
equations (SNSE) with additive noise for velocity $v$ and pressure $p$ in a
viscous incompressible flow:%
\begin{eqnarray}
dv(t) &=&\left[ \frac{\sigma ^{2}}{2}\Delta v-(v,\nabla )v-\nabla p+f(t,x)%
\right] dt+\sum_{r=1}^{q}\gamma _{r}(t,x)dw_{r}(t),  \label{sns1} \\
\ \ 0 &<&t\leq T,\ x\in \mathbf{R}^{2},  \notag \\
\mathop{\rm div}v &=&0,  \label{sns3}
\end{eqnarray}%
with spatial periodic conditions%
\begin{eqnarray}
v(t,x+Le_{i}) &=&v(t,x),\ p(t,x+Le_{i})=p(t,x),  \label{sns4} \\
0 &\leq &t\leq T,\ \ i=1,2,  \notag
\end{eqnarray}%
and the initial condition%
\begin{equation}
v(0,x)=\varphi (x).  \label{sns5}
\end{equation}%
In (\ref{sns1})-(\ref{sns4}), $v,$ $f,$ and $\gamma _{r}$ are
two-dimensional functions;$\ p$ is a scalar; $\{e_{i}\}$ is the canonical
basis in $\mathbf{R}^{2}$ and $L>0$ is the period. The functions $f=f(t,x)$
and $\gamma _{r}(t,x)$ are assumed to be spatial periodic as well. Further,
we require that $\gamma _{r}(t,x)$ are divergence free:
\begin{equation}
\mathop{\rm div}\gamma _{r}(t,x)=0,\ r=1,\ldots ,q.  \label{sns6}
\end{equation}%
For simplicity of proofs, we assume that the number of noises $q$ is finite
but it can be shown that the theoretical results of this section are also
valid when $q$ is infinite if $\Vert \gamma _{r}(t,x)\Vert _{m}$ for some $%
m\geq 0$ decay exponentially fast with increase of $r.$

\medskip

\noindent \textbf{Assumption 5.1.} \textit{We assume that the coefficients} $%
f(t,x)$ \textit{and} $\gamma _{r}(s,x),$ $r=1,\ldots ,q,$ \textit{belong to}
$\mathbf{H}_{p}^{m+1}(Q)$ \textit{and the initial condition} $\varphi (x)$
\textit{belongs to} $\mathbf{H}_{p}^{m+2}(Q)$ \textit{for some} $m\geq 0$.

\medskip

Under this assumption the problem (\ref{sns1})-(\ref{sns5})\ has a unique
solution $v(t,x),\ p(t,x),$ $(t,x)\in \lbrack 0,T]\times R^{2},$ so that for
some $m\geq 0$ and $l\geq 2$ \cite{MatPhD,Mat02}:%
\begin{equation}
E\Vert v(t,\cdot )\Vert _{m+2}^{l}\leq K,  \label{momv}
\end{equation}%
where $K>0$ may depend on $l$, $m,\ T,$ $f(t,x),$ $\gamma _{r}(t,x),$ and $%
\varphi (x).$ The solution $v(t,x),\ p(t,x),$ $(t,x)\in \lbrack 0,T]\times
\mathbf{R}^{2},$ to (\ref{sns1})-(\ref{sns5}) is $\mathcal{F}_{t}$-adaptive,
$v(t,\cdot )\in \mathbf{V}_{p}^{m+2}$ and $\nabla p(t,\cdot )\in (\mathbf{V}%
_{p}^{m+2})^{\bot }$ for every $t\in \lbrack 0,T]$ and $\mathbf{\omega }\in
\Omega .$ We note that if we were interested in variational solutions of (%
\ref{sns1})-(\ref{sns5}) then it is more natural to put $m\geq -1$ in
Assumption~5.1; but here our focus is on the vorticity formulation and then
it is natural to require more, $m\geq 0$.

The vorticity formulation of the problem (\ref{sns1})-(\ref{sns5}) has the
form
\begin{equation}
d\omega =\left[ \frac{\sigma ^{2}}{2}\Delta \omega -(v,\nabla )\omega +g(t,x)%
\right] dt+\sum_{r=1}^{q}\mu _{r}(t,x)dw_{r}(t),  \label{sns7}
\end{equation}%
where $g=\mathop{\rm curl}f$ and $\mu _{r}=\mathop{\rm curl}\gamma _{r}.$
The vorticity satisfies the initial and periodic boundary conditions
\begin{equation}
\omega (0,x)=\mathop{\rm curl}\varphi (x):=\phi (x)  \label{sns8}
\end{equation}%
and spatial periodic conditions%
\begin{equation}
\omega (t,x+Le_{i})=\omega (t,x),\ i=1,2,\ 0\leq t\leq T.  \label{sns9}
\end{equation}%
We note that $\omega (t,x)$ is a one-dimensional function here. Using the
linear operator $U$ from (\ref{Mnew1}), we can re-write (\ref{sns7}) as
\begin{equation}
d\omega =\left[ \frac{\sigma ^{2}}{2}\Delta \omega -(U\omega ,\nabla )\omega
+g(t,x)\right] dt+\sum_{r=1}^{q}\mu _{r}(t,x)dw_{r}(t).  \label{sns10}
\end{equation}

Similarly to the solution $v(t,x)$ of (\ref{sns1})-(\ref{sns5}), the
solution $\omega (t,x)$ to the vorticity problem (\ref{sns7})-(\ref{sns9})
under Assumption~5.1 is so that for some $m\geq 0$ and $p\geq 2$:%
\begin{equation}
E\Vert \omega (t,\cdot )\Vert _{m+1}^{p}\leq K,  \label{sns100}
\end{equation}%
where $K>0$ depends on $p$, $m,$ $g,$ $\mu _{r},$ and $\phi .$ Note that
under Assumption~5.1 the coefficients $g(t,x)$ and $\mu _{r}(s,x),$ $%
r=1,\ldots ,q,$ belong to $\mathbf{H}_{p}^{m}(Q)$ and the initial condition $%
\phi (x)$ belongs to $\mathbf{H}_{p}^{m+1}(Q).$ As it is clear from the
context, we are dealing here with solutions understood in the strong sense
probabilistically and PDE-wise in the variational sense.

\subsection{Two technical lemmas\label{sec:sns1}}

For proving convergence of the numerical method in Section~\ref{sec:sns4},
we need two further properties of the solution $\omega (t,x)$ which are
formulated in the next two lemmas.

It is convenient to introduce the notation for the solution $\omega (t,x)$
of the problem (\ref{sns7})-(\ref{sns9}) which reflects its dependence on
the initial condition $\phi (x)$ prescribed at time $s\leq t$:
\begin{equation*}
\omega (t,x)=\omega (t,x;s,\phi ).
\end{equation*}

Let us prove a technical lemma which is related to Lemmas~4.10(1) and A.1
from \cite{HaMa06}.

\begin{lemma}
\label{LemMat}Let Assumption~5.1 hold with $m=0.$ There exist constants $%
\beta _{0}>0$ and $\alpha >0$ such that for any $\beta \in (0,\beta _{0}]$
and $0\leq t\leq t+h\leq T:$
\begin{eqnarray}
&&E\exp \left( \beta \left[ ||\omega (t+h,\cdot ;t,\phi )||^{2}-||\phi ||^{2}%
\right] +\beta \frac{\sigma ^{2}}{4}\int_{t}^{t+h}||\nabla \omega (s,\cdot
;t,\phi )||^{2}ds\right)  \label{sns30x} \\
&\leq &\exp \left( \beta \int_{t}^{t+h}\left( \frac{2}{\alpha \sigma ^{2}}%
||g(s,\cdot )||^{2}+\sum_{r=1}^{q}||\mu _{r}(s,\cdot )||^{2}\right)
ds\right) .  \notag
\end{eqnarray}
\end{lemma}

\noindent \textbf{Proof}. By the Ito formula, integration by parts and using
$\mathop{\rm div}v(t,x)=0$, we obtain
\begin{gather}
\frac{1}{2}d||\omega (s,\cdot )||^{2}=\left[ -\frac{\sigma ^{2}}{2}||\nabla
\omega (s,\cdot )||^{2}+(g(s,\cdot ),\omega (s,\cdot ))+\frac{1}{2}%
\sum_{r=1}^{q}||\mu _{r}(s,\cdot )||^{2}\right] ds  \label{sns302} \\
+\sum_{r=1}^{q}(\mu _{r}(s,\cdot ),\omega (s,\cdot ))dw_{r}(s),\ t<s\leq t+h,
\notag \\
||\omega (t,\cdot )||^{2}=||\phi ||^{2}.  \notag
\end{gather}%
Using the elementary inequality, we get for any $\alpha >0:$
\begin{eqnarray}
&&\frac{1}{2}d||\omega (s,\cdot )||^{2}  \label{sns304} \\
&\leq &\left[ -\frac{\sigma ^{2}}{2}||\nabla \omega (s,\cdot )||^{2}+\frac{1%
}{\alpha \sigma ^{2}}||g(s,\cdot )||^{2}+\frac{\sigma ^{2}}{4}\alpha
||\omega (s,\cdot )||^{2}+\frac{1}{2}\sum_{r=1}^{q}||\mu _{r}(s,\cdot )||^{2}%
\right] ds  \notag \\
&&+\sum_{r=1}^{q}(\mu _{r}(s,\cdot ),\omega (s,\cdot ))dw_{r}(s).  \notag
\end{eqnarray}%
By Poincaire's inequality (\ref{Poin}), for some $\alpha >0,$ we have
\begin{equation}
||\nabla \omega (t,\cdot )||^{2}\geq \alpha ||\omega (t,\cdot )||^{2}.
\label{sns303}
\end{equation}%
By (\ref{sns303}), we obtain
\begin{eqnarray}
d||\omega (s,\cdot )||^{2} &\leq &\left[ -\frac{\sigma ^{2}}{4}||\nabla
\omega (s,\cdot )||^{2}-\frac{\sigma ^{2}}{4}\alpha ||\omega (s,\cdot
)||^{2}+\frac{2}{\alpha \sigma ^{2}}||g(s,\cdot )||^{2}\right.
\label{sns3045} \\
&&\left. +\sum_{r=1}^{q}||\mu _{r}(s,\cdot )||^{2}\right] ds+2%
\sum_{r=1}^{q}(\mu _{r}(s,\cdot ),\omega (s,\cdot ))dw_{r}(s),  \notag
\end{eqnarray}%
then for any $c>0$%
\begin{eqnarray}
&&c||\omega (t+h,\cdot )||^{2}-c||\phi ||^{2}+c\frac{\sigma ^{2}}{4}%
\int_{t}^{t+h}||\nabla \omega (s,\cdot )||^{2}ds  \label{sns306} \\
&&-c\int_{t}^{t+h}\left( \frac{2}{\alpha \sigma ^{2}}||g(s,\cdot
)||^{2}+\sum_{r=1}^{q}||\mu _{r}(s,\cdot )||^{2}\right) ds  \notag \\
&\leq &2c\int_{t}^{t+h}\sum_{r=1}^{q}(\mu _{r}(s,\cdot ),\omega (s,\cdot
))dw_{r}(s)-\alpha \frac{\sigma ^{2}}{4}c\int_{t}^{t+h}||\omega (s,\cdot
)||^{2}ds.  \notag
\end{eqnarray}%
Let
\begin{equation*}
M(t,t^{\prime }):=2c\int_{t}^{t^{\prime }}\sum_{r=1}^{q}(\mu _{r}(s,\cdot
),\omega (s,\cdot ))dw_{r}(s)
\end{equation*}%
which is a continuous $L^{2}$-martingale with quadratic variation
\begin{equation*}
<M>(t,t^{\prime }):=4c^{2}\int_{t}^{t^{\prime }}\sum_{r=1}^{q}(\mu
_{r}(s,\cdot ),\omega (s,\cdot ))^{2}ds.
\end{equation*}%
There exists a constant $\beta _{0}>0$ (independent of $h$ and $c)$ so that
for all $\beta \in (0,\beta _{0}]:$%
\begin{equation*}
\alpha \frac{\sigma ^{2}}{4}c\int_{t}^{t^{\prime }}||\omega (s,\cdot
)||^{2}ds\geq \frac{\beta }{2c}<M>(t,t^{\prime }).
\end{equation*}%
Hence
\begin{eqnarray}
&&c||\omega (t+h,\cdot )||^{2}-c||\phi ||^{2}+c\frac{\sigma ^{2}}{4}%
\int_{t}^{t+h}||\nabla \omega (s,\cdot )||^{2}ds  \label{sns305} \\
&&-c\int_{t}^{t+h}\left( \frac{2}{\alpha \sigma ^{2}}||g(s,\cdot
)||^{2}+\sum_{r=1}^{q}||\mu _{r}(s,\cdot )||^{2}\right) ds  \notag \\
&\leq &M(t,t+h)-\frac{\beta }{2c}<M>(t,t+h).  \notag
\end{eqnarray}%
For $c=\beta ,$ the right-hand side of (\ref{sns305}) is logarithm of a
local exponential martingale and therefore
\begin{gather*}
E\exp \left[ \beta ||\omega (t+h\wedge \tau _{n},\cdot )||^{2}-\beta ||\phi
||^{2}+\beta \frac{\sigma ^{2}}{4}\int_{t}^{t+h\wedge \tau _{n}}||\nabla
\omega (s,\cdot )||^{2}ds\right.  \\
\left. -\beta \int_{t}^{t+h\wedge \tau _{n}}\left( \frac{2}{\alpha \sigma
^{2}}||g(s,\cdot )||^{2}+\sum_{r=1}^{q}||\mu _{r}(s,\cdot )||^{2}\right) ds%
\right] \leq 1,
\end{gather*}%
where $\tau _{n}=\inf \{s>t:<M>(t,s)\geq n\}$ for a natural number $n.$
Tending $n$ to infinity, we arrive at (\ref{sns30x}). Lemma~\ref{LemMat} is
proved. \bigskip

Note that it follows from (\ref{sns30x}) that
\begin{eqnarray}
&&E\exp \left( \beta \frac{\sigma ^{2}}{4}\int_{t}^{t+h}||\nabla \omega
(s,\cdot ;t,\phi )||^{2}ds\right)  \label{sns30xx} \\
&\leq &\exp \left( \beta ||\phi ||^{2}+\beta \int_{t}^{t+h}\left( \frac{2}{%
\alpha \sigma ^{2}}||g(s,\cdot )||^{2}+\sum_{r=1}^{q}||\mu _{r}(s,\cdot
)||^{2}\right) ds\right) .  \notag
\end{eqnarray}%
We also pay attention that the prove of Lemma~\ref{LemMat} is not relying on
smallness of the time step $h$ and, after replacing $t$ with $0$ and $t+h$
with $T,$ the result remains valid:
\begin{eqnarray}
&&E\exp \left( \beta \frac{\sigma ^{2}}{4}\int_{0}^{T}||\nabla \omega
(s,\cdot ;0,\phi )||^{2}ds\right)  \label{sns30xxg} \\
&\leq &\exp \left( \beta ||\phi ||^{2}+\beta \int_{0}^{T}\left( \frac{2}{%
\alpha \sigma ^{2}}||g(s,\cdot )||^{2}+\sum_{r=1}^{q}||\mu _{r}(s,\cdot
)||^{2}\right) ds\right) .  \notag
\end{eqnarray}

We now prove the next lemma which gives us dependence of the solution $%
\omega (s,x;t,\phi )$ on the initial data.

\begin{lemma}
\label{Lem2}Let Assumption~5.1 hold with $m=2$ and $\phi _{i}(t,x),$ $i=1,2,$
be $\mathcal{F}_{t}$-measurable processes satisfying (\ref{sns100}) with $%
m=2 $. There exists a constant $c_{0}>0$ such that for every $c\in (0,c_{0})$
there is a sufficiently small $h>0$ so that we have for $t\leq s\leq t+h:$%
\begin{equation}
\omega (s,x;t,\phi _{1}(t,\cdot ))-\omega (s,x;t,\phi _{2}(t,\cdot ))=\phi
_{1}(t,x)-\phi _{2}(t,x)+\eta (s,x)  \label{sns30}
\end{equation}%
for which
\begin{eqnarray}
&&||\omega (s,\cdot ;t,\phi _{1})-\omega (s,\cdot ;t,\phi _{2})||^{2}
\label{sns31} \\
&\leq &||\phi _{1}(t,\cdot )-\phi _{2}(t,\cdot )||^{2}\exp \left(
K(s-t)+c\int_{t}^{s}||\nabla \omega (s^{\prime },\cdot ;t,\phi _{1}(t,\cdot
))||^{2}ds^{\prime }\right) \mathbf{,}  \notag
\end{eqnarray}%
where $K>0$ is a constant.

The process $\eta (s)$ satisfies the following estimate
\begin{equation}
||\eta (s,\cdot )||^{2}\leq (s-t)||\phi _{1}(t,\cdot )-\phi _{2}(t,\cdot
)||^{2}+C(s,\mathbf{\omega })(s-t)^{3},  \label{sns32}
\end{equation}%
where $C(s,\mathbf{\omega })>0$ is an $\mathcal{F}_{s}$-adapted process with
bounded moments of a sufficiently high order.
\end{lemma}

\noindent \textbf{Proof}. Let
\begin{equation*}
\theta (s,x):=\omega (s,x;t,\phi _{1})-\omega (s,x;t,\phi _{2})
\end{equation*}%
We have%
\begin{eqnarray*}
d\theta (s,x) &=&\left[ \frac{\sigma ^{2}}{2}\Delta \theta -(U\theta ,\nabla
)\omega (s,\cdot ;t,\phi _{1})-(U\omega (s,\cdot ;t,\phi _{2}),\nabla
)\theta \right] ds,\ t<s\leq t+h, \\
\theta (t,x) &=&\phi _{1}(t,x)-\phi _{2}(t,x).
\end{eqnarray*}%
Then
\begin{gather}
\frac{1}{2}d||\theta (s,\cdot )||^{2}=\left[ -\frac{\sigma ^{2}}{2}||\nabla
\theta (s,\cdot )||^{2}-((U\theta (s,\cdot ),\nabla )\omega (s,\cdot ;t,\phi
_{1}),\theta (s,\cdot ))\right] ds,\ t<s\leq t+h,\ \   \label{sns33} \\
||\theta (t,\cdot )||^{2}=||\phi _{1}(t,\cdot )-\phi _{2}(t,\cdot )||^{2}.
\notag
\end{gather}%
Using the inequality (\ref{A66}) with $c_{2}=\sigma ^{2}/4$, we have that
there exists $K>0$ such that for any $c>0$
\begin{eqnarray}
2|((U\theta (s,\cdot ),\nabla )\omega (s,\cdot ;t,\phi _{1}),\theta (s,\cdot
))| &\leq &\frac{\sigma ^{2}}{2}||\nabla \theta (s,\cdot )||^{2}+K||\theta
(s,\cdot )||^{2}  \label{sns34} \\
&&+c||\nabla \omega (s,\cdot ;t,\phi _{1})||^{2}||\theta (s,\cdot )||^{2}
\notag
\end{eqnarray}%
and hence%
\begin{equation*}
d||\theta (s,\cdot )||^{2}\leq \left[ K+c||\nabla \omega (s,\cdot ;t,\phi
_{1})||^{2}\right] ||\theta (s,\cdot )||^{2}ds,\ t<s\leq t+h,
\end{equation*}%
which implies
\begin{equation}
||\theta (s,\cdot )||^{2}\leq ||\phi _{1}(t,\cdot )-\phi _{2}(t,\cdot
)||^{2}\exp \left( K(s-t)+c\int_{t}^{s}||\nabla \omega (s^{\prime },\cdot
;t,\phi _{1}(t,\cdot ))||^{2}ds^{\prime }\right) .  \label{sns35}
\end{equation}%
Thus we have proved the inequality (\ref{sns31}).

Let us now prove (\ref{sns32}). We have
\begin{equation*}
\eta (s,x)=\int_{t}^{s}\left[ \frac{\sigma ^{2}}{2}\Delta \theta -(U\theta
,\nabla )\omega (s^{\prime },\cdot ;t,\phi _{1})-(U\omega (s^{\prime },\cdot
;t,\phi _{2}),\nabla )\theta \right] ds^{\prime },\ t<s\leq t+h,
\end{equation*}%
which together with (\ref{momv}) and (\ref{sns100}) implies that
\begin{equation}
||\eta (s,x)||\leq C(s,\mathbf{\omega })(s-t),  \label{sns37}
\end{equation}%
where $C(s,\mathbf{\omega })>0$ is an $\mathcal{F}_{s}$-adapted process with
bounded moments of a sufficiently high order. It is not difficult to see
that the inequality (\ref{sns37}) is also valid for $||\nabla \eta (s,x)||$
and $||\Delta \eta (s,x)||:$%
\begin{equation}
||\nabla \eta (s,x)||\leq C(s,\mathbf{\omega })(s-t),\ \ ||\Delta \eta
(s,x)||\leq C(s,\mathbf{\omega })(s-t).  \label{sns377}
\end{equation}%
We have
\begin{equation*}
d||\eta (s^{\prime },x)||^{2}=\left[ (\sigma ^{2}\Delta \theta ,\eta
)-2((U\theta ,\nabla )\omega (s^{\prime },\cdot ;t,\phi _{1}),\eta
)-2((U\omega (s^{\prime },\cdot ;t,\phi _{2}),\nabla )\theta ,\eta )\right]
ds^{\prime }.
\end{equation*}%
Using integration by parts, (\ref{sns31}), and (\ref{sns377}) (we also
recall that $s^{\prime }-t\leq h$ which is sufficiently small), we get
\begin{equation*}
|(\sigma ^{2}\Delta \theta ,\eta )|=\sigma ^{2}|(\theta ,\Delta \eta )|\leq
\sigma ^{2}||\theta ||||\Delta \eta ||\leq C(s^{\prime },\mathbf{\omega }%
)(s^{\prime }-t)||\phi _{1}(t,\cdot )-\phi _{2}(t,\cdot )||.
\end{equation*}%
By (\ref{A2}) with $m_{1}=1,$ $m_{2}=0$, and $m_{3}=1$, (\ref{A3}), (\ref%
{sns100}), (\ref{sns31}), (\ref{sns37}) and (\ref{sns377})$,$ we obtain
\begin{eqnarray*}
|2((U\theta ,\nabla )\omega (s,\cdot ;t,\phi _{1}),\eta )| &\leq &K||U\theta
||_{1}||\omega ||_{1}||\eta ||_{1}\leq K||\theta ||||\omega ||_{1}||\eta
||_{1} \\
&\leq &C(s^{\prime },\mathbf{\omega })(s^{\prime }-t)||\phi _{1}(t,\cdot
)-\phi _{2}(t,\cdot )||.
\end{eqnarray*}%
And by (\ref{A2}) with $m_{1}=1,$ $m_{2}=1$, and $m_{3}=0$, (\ref{A3}), (\ref%
{sns100}), (\ref{sns31}), (\ref{sns37}) and (\ref{sns377}), we arrive at
\begin{eqnarray*}
|2((U\omega (s^{\prime },\cdot ;t,\phi _{2}),\nabla )\theta ,\eta )|
&=&2|((U\omega (s^{\prime },\cdot ;t,\phi _{2}),\nabla )\eta ,\theta )|\leq
K||\omega ||||\eta ||_{2}||\theta || \\
&\leq &C(s^{\prime },\mathbf{\omega })(s^{\prime }-t)||\phi _{1}(t,\cdot
)-\phi _{2}(t,\cdot )||.
\end{eqnarray*}%
Then we have
\begin{eqnarray*}
d||\eta (s^{\prime },x)||^{2} &\leq &C(s^{\prime },\mathbf{\omega }%
)(s^{\prime }-t)||\phi _{1}(t,\cdot )-\phi _{2}(t,\cdot )||ds^{\prime } \\
&\leq &||\phi _{1}(t,\cdot )-\phi _{2}(t,\cdot )||^{2}ds^{\prime }+\frac{%
C^{2}(s^{\prime },\mathbf{\omega })}{4}(s^{\prime }-t)^{2}ds^{\prime }
\end{eqnarray*}%
from which (\ref{sns32}) follows. \ Lemma~\ref{Lem2} is proved.

\subsection{One-step approximation\label{sec:sns2}}

Similarly to derivation of the approximation for the deterministic NSE in
Section~\ref{sec:approx}, we can approximate the stochastic NSE (\ref{sns7}%
)-(\ref{sns9}) by freezing the velocity as in (\ref{eq:freeze}):
\begin{equation}
v(t,x)\approx \hat{v}(t,x):=v(t_{k},x):=\hat{v}(x),\ t_{k}<t\leq t_{k+1},
\label{eq:freeze2}
\end{equation}%
and obtain an approximation $\tilde{\omega}(t,x)$ of $\omega (t,x)$ on $\
t_{k}\leq t\leq t_{k+1},$ as follows
\begin{eqnarray}
d\tilde{\omega}=\left[ \frac{\sigma ^{2}}{2}\Delta \tilde{\omega}-(\hat{v}%
,\nabla )\tilde{\omega}+g(t,x)\right] dt+\sum_{r=1}^{q}\mu
_{r}(t,x)dw_{r}(t),\ t_{k}<t\leq t_{k+1},  \label{sns11} \\
\tilde{\omega}(t_{k},x)=\omega (t_{k},x),\ \ \tilde{\omega}(t_{k},x+Le_{j})=%
\tilde{\omega}(t_{k},x),\ j=1,2.  \label{sns12}
\end{eqnarray}

It is not difficult to see that the local error $\delta _{\omega }(t,x)=%
\tilde{\omega}(t,x)-\omega (t,x),$ $t_{k}\leq t\leq t_{k+1},$ for the
approximation $\tilde{\omega}(t,x)$ of the solution $\omega (t,x)$ of the
stochastic NSE (\ref{sns7})-(\ref{sns9}) satisfies the problem of the same
form as (\ref{M8})-(\ref{M9}) but with positive direction of time:
\begin{eqnarray}
d\delta _{\omega }=\left[ \frac{\sigma ^{2}}{2}\Delta \delta _{\omega
}-(v,\nabla )\delta _{\omega }-((v-\hat{v}),\nabla )\tilde{\omega}\right] dt,
\label{sns14} \\
\delta _{\omega }(t_{k},x)=0.  \label{sns15}
\end{eqnarray}

We note that the main difference of (\ref{sns14})-(\ref{sns15}) with (\ref%
{M8})-(\ref{M9}) is that the functions in (\ref{sns14}) are random and
non-smooth in time, they have the same regularity in time as Wiener
processes.

Moments of $||\tilde{\omega}||_{3}$ (and hence of $||\delta _{\omega
}||_{3}) $ up to a sufficiently high order are bounded under Assumption~5.1
with $m=2$: for $t_{k}<t\leq t_{k+1}$ and $p\geq 1:$%
\begin{equation}
E\Vert \tilde{\omega}(t,\cdot )\Vert _{3}^{2p}\leq K,  \label{sns159}
\end{equation}%
where $K>0$ is a constant, which can be proved by arguments similar to
boundedness of the global approximation (see Theorems~\ref{thm:sns2} and~\ref%
{thm:sns2nn}) but not considered here.

To obtain bounds for the one-step error $\delta _{\omega },$ we first prove
the following lemma.

\begin{lemma}
\label{Lem1} Let Assumption~5.1 hold with $m=1.$ For $v(t,x)$ from $(\ref%
{sns1})$-$(\ref{sns5})$, $\hat{v}(x)$ from $(\ref{eq:freeze2})$, and $\tilde{%
\omega}(t,x)$ from $(\ref{sns11})$-$(\ref{sns12})$, we have for $t_{k}<t\leq
t_{k+1}$ and sufficiently small $h>0:$
\begin{eqnarray}
||E[((v-\hat{v}),\nabla )\tilde{\omega}|\mathcal{F}_{t_{k}}]|| &\leq
&C(t_{k},\mathbf{\omega })h,  \label{sns16} \\
\left( E||v-\hat{v}||^{2p}\right) ^{1/2p} &\leq &Kh^{1/2},\ p\geq 1,
\label{sns17}
\end{eqnarray}%
where $C(t_{k},\mathbf{\omega })>0$ is an $\mathcal{F}_{t_{k}}$-measurable
random variable with moments of a sufficiently high order bounded by a
constant independent of $h$ and $K>0$ is a constant independent of $h.$
\end{lemma}

\noindent \textbf{Proof}. \ From (\ref{sns1}) and (\ref{eq:freeze2}), we
have for $t_{k}<t\leq t_{k+1}:$%
\begin{equation}
v(t,x)-\hat{v}(x)=\int_{t_{k}}^{t}\left[ \frac{\sigma ^{2}}{2}\Delta
v-(v,\nabla )v-\nabla p+f(s,x)\right] ds+\int_{t_{k}}^{t}\sum_{r=1}^{q}%
\gamma _{r}(s,x)dw_{r}(s).  \label{sns177}
\end{equation}%
Then it is not difficult to obtain the estimate (\ref{sns17}) using (\ref%
{momv}) and the assumptions on $f$ and $\gamma _{r}.$

From (\ref{sns177}) and (\ref{sns11}), we have
\begin{eqnarray*}
((v-\hat{v}),\nabla )\tilde{\omega} &=&\left( \int_{t_{k}}^{t}\left[ \frac{%
\sigma ^{2}}{2}\Delta v-(v,\nabla )v-\nabla p+f(s,x)\right] ds,\nabla
\right) \tilde{\omega}(t,x) \\
&&+\left( \int_{t_{k}}^{t}\sum_{r=1}^{q}\gamma _{r}(s,x)dw_{r}(s),\nabla
\right) \\
&&\left\{ \tilde{\omega}(t_{k},x)+\int_{t_{k}}^{t}\left[ \frac{\sigma ^{2}}{2%
}\Delta \tilde{\omega}-(\hat{v},\nabla )\tilde{\omega}+g(s,x)\right]
ds+\int_{t_{k}}^{t}\sum_{r=1}^{q}\mu _{r}(s,x)dw_{r}(s)\right\}
\end{eqnarray*}%
from which it is not difficult to see that the inequality (\ref{sns16})
holds. Lemma~\ref{Lem1} is proved. \medskip

Now we proceed to proving estimates for the one-step error of $\tilde{\omega}%
(t,x).$

\begin{theorem}
\label{thm:sns1}Let Assumption~5.1 hold with $m=2$. The one-step error of $%
\tilde{\omega}(t,x),\ t_{k}\leq t\leq t_{k+1},$\ which solves $(\ref{sns11})$%
-$(\ref{sns12}),$ has the following bounds for $t_{k}\leq t\leq t_{k+1}$ and
sufficiently small $h>0:$%
\begin{eqnarray}
||E[\delta _{\omega }(t,\cdot )|\mathcal{F}_{t_{k}}]|| &\leq &C(t_{k},%
\mathbf{\omega })h^{2},\   \label{sns18} \\
\left( E||\delta _{\omega }(t,\cdot )||^{2}\right) ^{1/2} &\leq &Kh^{3/2},\
\label{sns19}
\end{eqnarray}%
where $C(t_{k},\mathbf{\omega })>0$ is an $\mathcal{F}_{t_{k}}$-measurable
random variable with moments of a sufficiently high order bounded by a
constant independent of $h$ and $K>0$ is a constant independent of $h.$
\end{theorem}

\noindent \textbf{Proof}. Taking scalar product of (\ref{sns14}) and $\delta
_{\omega }(t,x),$ using integration by parts and the property $%
\mathop{\rm
div}v(t,x)=0$, we get
\begin{eqnarray}
\frac{1}{2}d||\delta _{\omega }(t,\cdot )||^{2} &=&\frac{\sigma ^{2}}{2}%
(\Delta \delta _{\omega }(t,\cdot ),\delta _{\omega }(t,\cdot ))dt-(\left[
(v(t,\cdot ),\nabla )\delta _{\omega }(t,\cdot )\right] ,\delta _{\omega
}(t,\cdot ))dt  \label{sns200} \\
&&-(\left[ ((v(t,\cdot )-\hat{v}(\cdot )),\nabla )\tilde{\omega}(t,\cdot )%
\right] ,\delta _{\omega }(t,\cdot ))dt  \notag \\
&=&-\frac{\sigma ^{2}}{2}||\nabla \delta _{\omega }(t,\cdot )||^{2}dt-(\left[
((v(t,\cdot )-\hat{v}(\cdot )),\nabla )\tilde{\omega}(t,\cdot )\right]
,\delta _{\omega }(t,\cdot ))dt.  \notag
\end{eqnarray}%
Then
\begin{equation}
\frac{1}{2}dE||\delta _{\omega }(t,\cdot )||^{2}=-\frac{\sigma ^{2}}{2}%
E||\nabla \delta _{\omega }(t,\cdot )||^{2}dt-E(\left[ ((v(t,\cdot )-\hat{v}%
(\cdot )),\nabla )\tilde{\omega}(t,\cdot )\right] ,\delta _{\omega }(t,\cdot
))dt.  \label{sns20}
\end{equation}%
For the last term in (\ref{sns20}), we get
\begin{eqnarray}
&&|E(\left[ ((v(t,\cdot )-\hat{v}(\cdot )),\nabla )\tilde{\omega}(t,\cdot )%
\right] ,\delta _{\omega }(t,\cdot ))|  \label{sns201} \\
&\leq &KE||v(t,\cdot )-\hat{v}(\cdot )||\cdot ||\tilde{\omega}(t,\cdot
)||_{3}\cdot ||\delta _{\omega }(t,\cdot )||  \notag \\
&\leq &K\left( E||v(t,\cdot )-\hat{v}(\cdot )||^{2}\cdot ||\tilde{\omega}%
(t,\cdot )||_{3}^{2}\right) ^{1/2}\left( E||\delta _{\omega }(t,\cdot
)||^{2}\right) ^{1/2}  \notag \\
&\leq &K\left( E||v(t,\cdot )-\hat{v}(\cdot )||^{4}\right) ^{1/4}\left( E||%
\tilde{\omega}(t,\cdot )||_{3}^{4}\right) ^{1/4}\left( E||\delta _{\omega
}(t,\cdot )||^{2}\right) ^{1/2}  \notag \\
&\leq &Kh^{1/2}\left( E||\delta _{\omega }(t,\cdot )||^{2}\right) ^{1/2},
\notag
\end{eqnarray}%
where for the first line we used the inequality (\ref{A2}) with $m_{1}=0,$ $%
m_{2}=2,$ $m_{3}=0$; we applied the Cauchy-Bunyakovski inequality twice to
arrive at the pre-last line; and we used the error estimate (\ref{sns17})
with $p=2$ and boundedness of the moment $E||\tilde{\omega}(t,\cdot
)||_{3}^{4}$ (see (\ref{sns159})) to obtain the last line.

Thus
\begin{equation*}
\frac{1}{2}dE||\delta _{\omega }(t,\cdot )||^{2}\leq Kh^{1/2}\left(
E||\delta _{\omega }(t,\cdot )||^{2}\right) ^{1/2}dt,
\end{equation*}%
and since $\delta _{\omega }(t_{k},x)=0,$ we arrive at
\begin{equation*}
\int_{t_{k}}^{t}\frac{1}{2}\frac{dE||\delta _{\omega }(s,\cdot )||^{2}}{%
\left( E||\delta _{\omega }(s,\cdot )||^{2}\right) ^{1/2}}=\left( E||\delta
_{\omega }(t,\cdot )||^{2}\right) ^{1/2}\leq Kh^{3/2}
\end{equation*}%
confirming (\ref{sns19}).

Now we are to prove (\ref{sns18}). Using (\ref{sns14}), we write the
equation for $dE[\delta _{\omega }(t,x)|\mathcal{F}_{t_{k}}]$ and, after
taking scalar product of the components of this equation and $E[\delta
_{\omega }(t,x)|\mathcal{F}_{t_{k}}]$ and doing integration by parts, we
arrive
\begin{eqnarray}
&&\frac{1}{2}d||E[\delta _{\omega }(t,\cdot )|\mathcal{F}_{t_{k}}]||^{2}
\label{sns21} \\
&=&-\frac{\sigma ^{2}}{2}||\nabla E[\delta _{\omega }(t,\cdot )|\mathcal{F}%
_{t_{k}}]||^{2}dt-(E\left[ (v(t,\cdot ),\nabla )\delta _{\omega }(t,\cdot )|%
\mathcal{F}_{t_{k}}\right] ,E[\delta _{\omega }(t,\cdot )|\mathcal{F}%
_{t_{k}}])dt\   \notag \\
&&-(E\left[ ((v(t,\cdot )-\hat{v}(\cdot )),\nabla )\tilde{\omega}(t,\cdot )|%
\mathcal{F}_{t_{k}}\right] ,E[\delta _{\omega }(t,\cdot )|\mathcal{F}%
_{t_{k}}])dt.  \notag
\end{eqnarray}%
By (\ref{sns16}), we get for the third term in (\ref{sns21}):
\begin{eqnarray}
&&|(E\left[ ((v(t,\cdot )-\hat{v}(\cdot )),\nabla )\tilde{\omega}(t,\cdot )|%
\mathcal{F}_{t_{k}}\right] ,E[\delta _{\omega }(t,\cdot )|\mathcal{F}%
_{t_{k}}])|  \label{sns22} \\
&\leq &||E\left[ ((v(t,\cdot )-\hat{v}(\cdot )),\nabla )\tilde{\omega}%
(t,\cdot )|\mathcal{F}_{t_{k}}\right] ||\cdot ||E[\delta _{\omega }(t,\cdot
)|\mathcal{F}_{t_{k}}]||  \notag \\
&\leq &C(t_{k},\mathbf{\omega })h||E[\delta _{\omega }(t,\cdot )|\mathcal{F}%
_{t_{k}}]||,  \notag
\end{eqnarray}%
where $C(t_{k},\mathbf{\omega })>0$ is an $\mathcal{F}_{t_{k}}$-measurable
random variable which has moments up to a sufficiently high order and does
not depend on $h.$

By simple maniplations and using (\ref{A2}) $m_{1}=2,$ $m_{2}=0,$ $m_{3}=0$
as well as (\ref{A32}), the Cauchy-Bunyakovski inequality, (\ref{sns100}), a
conditional version of (\ref{sns19}), and (\ref{sns159}), we obtain for the
second term in (\ref{sns21}):%
\begin{eqnarray}
&&|(E\left[ (v(t,\cdot ),\nabla )\delta _{\omega }(t,\cdot )|\mathcal{F}%
_{t_{k}}\right] ,E[\delta _{\omega }(t,\cdot )|\mathcal{F}_{t_{k}}])|
\label{sns23} \\
&=&|E\left[ ((v(t,\cdot ),\nabla )\delta _{\omega }(t,\cdot ),E[\delta
_{\omega }(t,\cdot )|\mathcal{F}_{t_{k}}])|\mathcal{F}_{t_{k}}\right] |
\notag \\
&\leq &E\left[ |((v(t,\cdot ),\nabla )\delta _{\omega }(t,\cdot ),E[\delta
_{\omega }(t,\cdot )|\mathcal{F}_{t_{k}}])|\ |\mathcal{F}_{t_{k}}\right]
\notag \\
&=&E\left[ |((v(t,\cdot ),\nabla )E[\delta _{\omega }(t,\cdot )|\mathcal{F}%
_{t_{k}}],\delta _{\omega }(t,\cdot ))|\ |\mathcal{F}_{t_{k}}\right]   \notag
\\
&\leq &KE\left[ ||v(t,\cdot )||_{2}\cdot ||E[\delta _{\omega }(t,\cdot )|%
\mathcal{F}_{t_{k}}]||_{1}\cdot ||[\delta _{\omega }(t,\cdot )||\ |\mathcal{F%
}_{t_{k}}\right]   \notag \\
&=&K||E[\delta _{\omega }(t,\cdot )|\mathcal{F}_{t_{k}}]||_{1}E\left[
||v(t,\cdot )||_{2}\cdot ||[\delta _{\omega }(t,\cdot )||\ |\mathcal{F}%
_{t_{k}}\right]   \notag \\
&\leq &K||\nabla E[\delta _{\omega }(t,\cdot )|\mathcal{F}_{t_{k}}]||E\left[
||\omega (t,\cdot )||_{1}\cdot ||[\delta _{\omega }(t,\cdot )||\ |\mathcal{F}%
_{t_{k}}\right]   \notag \\
&\leq &K||\nabla E[\delta _{\omega }(t,\cdot )|\mathcal{F}_{t_{k}}]||(E\left[
||\omega (t,\cdot )||_{1}^{2}|\mathcal{F}_{t_{k}}\right] )^{1/2}(E\left[
||\delta _{\omega }(t,\cdot )||^{2}|\mathcal{F}_{t_{k}}\right] )^{1/2}
\notag \\
&\leq &C(t_{k},\mathbf{\omega })h^{3/2}||\nabla E[\delta _{\omega }(t,\cdot
)|\mathcal{F}_{t_{k}}]||,  \notag
\end{eqnarray}%
where $C(t_{k},\mathbf{\omega })>0$ is an $\mathcal{F}_{t_{k}}$-measurable
random variable which has moments up to a sufficiently high order and does
not depend on $h.$

Combining (\ref{sns21})-(\ref{sns23}), we arrive at
\begin{eqnarray*}
\frac{1}{2}d||E[\delta _{\omega }(t,\cdot )|\mathcal{F}_{t_{k}}]||^{2} &\leq
&-\frac{\sigma ^{2}}{2}||\nabla E[\delta _{\omega }(t,\cdot )|\mathcal{F}%
_{t_{k}}]||^{2}dt+ \\
&&+C(t_{k},\mathbf{\omega })h^{3/2}||\nabla E[\delta _{\omega }(t,\cdot )|%
\mathcal{F}_{t_{k}}]||dt+C(t_{k},\mathbf{\omega })h||E[\delta _{\omega
}(t,\cdot )|\mathcal{F}_{t_{k}}]||dt \\
&=&-\frac{1}{2}\left( \sigma ||\nabla E[\delta _{\omega }(t,\cdot )|\mathcal{%
F}_{t_{k}}]||-\frac{C(t_{k},\mathbf{\omega })}{\sigma }h^{3/2}\right) ^{2}dt+%
\frac{C^{2}(t_{k},\mathbf{\omega })}{2\sigma ^{2}}h^{3}dt \\
&&+C(t_{k},\mathbf{\omega })h||E[\delta _{\omega }(t,\cdot )|\mathcal{F}%
_{t_{k}}]||dt \\
&\leq &\frac{C^{2}(t_{k},\mathbf{\omega })}{2\sigma ^{2}}h^{3}dt+C(t_{k},%
\mathbf{\omega })h||E[\delta _{\omega }(t,\cdot )|\mathcal{F}_{t_{k}}]||dt.
\end{eqnarray*}%
Then, for some $\mathcal{F}_{t_{k}}$-measurable independent of $h$ random
variable $C(t_{k},\mathbf{\omega })>0,$ we have
\begin{equation*}
d||E[\delta _{\omega }(t,\cdot )|\mathcal{F}_{t_{k}}]||^{2}\leq C(t_{k},%
\mathbf{\omega })h^{3}dt+\frac{1}{h}||E[\delta _{\omega }(t,\cdot )|\mathcal{%
F}_{t_{k}}]||^{2}dt,
\end{equation*}%
from which (\ref{sns18}) follows taking into account that $||E[\delta
_{\omega }(t_{k},\cdot )|\mathcal{F}_{t_{k}}]||=0.$ Theorem~\ref{thm:sns1}
is proved. \medskip

As in the deterministic case, we define
\begin{equation}
\tilde{v}(t,x):=U\tilde{\omega}(t,x),\ \ t_{k}\leq t\leq t_{k+1},
\label{v_tilde}
\end{equation}%
where the operator $U$ is from (\ref{Mnew1}).

Using the idea of the proof of Corollary~\ref{Cor2D}, it is not difficult to
prove the following corollary to Theorem~\ref{thm:sns1}.

\begin{corollary}
\label{Corsns}The one-step error $\delta _{v}(t,x):=v(t,x)-\tilde{v}%
(t,x)=U\delta _{\omega }(t,x)$ of $v(t,x),\ t_{k}\leq t\leq t_{k+1},$\ has
the following bounds for $t_{k}\leq t\leq t_{k+1}:$%
\begin{eqnarray}
||E\delta _{v}(t,\cdot )|| &\leq &Kh^{2},\   \label{sns25} \\
\left( E||\delta _{v}(t,\cdot )||^{2}\right) ^{1/2} &\leq &Kh^{3/2},\
\label{sns26}
\end{eqnarray}%
where $K>0$ is independent of $h.$
\end{corollary}

\subsection{The method\label{sec:sns3}}

Analogously, how it was done in the deterministic case (see Section~\ref%
{sec:approxcon}), we can construct the global approximation for the
stochastic NSE (\ref{sns7})-(\ref{sns9}) based on the one-step approximation
(\ref{sns11})-(\ref{sns12}). On the first step of the method we set
\begin{equation*}
\tilde{\omega}(t_{0},x)=\mathop{\rm curl}v(t_{0},x)=\phi (x)=\mathop{\rm
curl}\varphi (x)
\end{equation*}%
and
\begin{equation*}
\hat{v}(x)=\hat{v}(t,x)=u(t_{0},x)=\varphi (x),\ \ 0=t_{0}\leq t\leq t_{1}.
\end{equation*}%
Then we solve the linear SPDE (\ref{sns11})-(\ref{sns12}) on\textit{\ }$%
[t_{0},t_{1}]$ to obtain $\tilde{\omega}(t,x)$ and to construct
\begin{equation*}
\hat{v}(t_{1},x)=U\tilde{\omega}(t_{1},x).
\end{equation*}%
On the second step we solve (\ref{sns11})-(\ref{sns12}) on\textit{\ }$%
[t_{1},t_{2}]$ having $\tilde{\omega}(t_{1},x)$ and setting $\hat{v}(t,x)=%
\hat{v}(x)=\hat{v}(t_{1},x)$ for $t_{1}<t\leq t_{2}.$ As a result, we obtain
$\tilde{\omega}(t,x)$ on $[t_{1},t_{2}]$ and $\hat{v}(t_{2},x)=U\tilde{\omega%
}(t_{2},x),$ and so on. Proceeding in this way, we obtain on the $N$-th step
the approximation $\tilde{\omega}(t,x)$ on $[t_{N-1},t_{N}]$ for $\omega
(t,x)$ having $\tilde{\omega}(t_{N-1},x)$ and $\hat{v}(x)=\hat{v}%
(t_{N-1},x)=U\tilde{\omega}(t_{N-1},x)$ and setting $\hat{v}(t,x)=\hat{v}(x)=%
\hat{v}(t_{N-1},x)$ for $t_{N-1}<t\leq t_{N}.$ Finally, $\hat{v}(T,x)=U%
\tilde{\omega}(T,x).$

In order to realise the approximation process described above, it is
sufficient that on every time interval $[t_{k},t_{k+1}],$ $k=0,\ldots ,N-1,$
there exists a solution of the linear SPDE (\ref{sns11})-(\ref{sns12}), we
denote such a solution $\tilde{\omega}_{k}(t,x)$ which satisfies the
condition
\begin{equation}
\tilde{\omega}_{k}(t_{k},x)=\left\{
\begin{array}{c}
\mathop{\rm curl}\varphi (x),\ k=0, \\
\tilde{\omega}_{k-1}(t_{k},x),\ k=1,\ldots ,N,%
\end{array}%
\right.  \label{omega_k}
\end{equation}%
and has the time-independent $\hat{u}(x)$ within each interval $%
(t_{k},t_{k+1}]$ $\ $defined as\
\begin{equation}
\hat{v}(x):=\hat{v}_{k}(x)=U\tilde{\omega}_{k}(t_{k},x),\ t_{k}<t\leq
t_{k+1}.  \label{v_k}
\end{equation}%
Clearly, $\hat{v}(x)$ used in (\ref{sns11}) are different on the time
intervals $(t_{k},t_{k+1}]$.

Before considering global errors of the approximation in Section~\ref%
{sec:sns4}, we now prove boundedness of the approximation's moments.

\begin{theorem}
\label{thm:sns2}Let Assumption~5.1 hold with $m=0$. The moments of the
global approximation $\tilde{\omega}_{k}(t_{k},x)$ and $\hat{v}_{k}(x)$ are
uniformly bounded in $h$ and $k:$
\begin{eqnarray}
E||\tilde{\omega}_{k}(t_{k+1},\cdot )||^{2p} &\leq &||\phi (\cdot
)||^{2p}+K\int_{0}^{t_{k+1}}\left( ||g(s,\cdot )||^{2p}+\sum_{r=1}^{q}||\mu
_{r}(s,\cdot )||^{2p}\right) ds,  \label{sns27} \\
E||\hat{v}_{k}(\cdot )||^{2p} &\leq &KE||\tilde{\omega}_{k}(t_{k},\cdot
)||^{2p},  \label{sns28}
\end{eqnarray}%
where $K>0$ is independent of $h$ and $t_{k}$ but depends on $p.$
\end{theorem}

\noindent \textbf{Proof}. \ For every sufficiently large integer $n$, define
the stopping time%
\begin{equation*}
\tau _{n}=\inf \{0<t\leq T:||\tilde{\omega}(t,\cdot )||^{2}\geq n\}.
\end{equation*}

Using the Ito formula, doing integration by parts and taking into account
that $\hat{v}_{k}(x)$ is divergence free, we obtain%
\begin{gather}
d||\tilde{\omega}_{k}(t,\cdot )||^{2p}=2p||\tilde{\omega}_{k}(t,\cdot
)||^{2(p-1)}  \label{sns292} \\
\cdot \left[ -\frac{\sigma ^{2}}{2}||\nabla \tilde{\omega}_{k}(t,\cdot
)||^{2}+(g(t,\cdot ),\tilde{\omega}_{k}(t,\cdot ))+\frac{2p-1}{2}%
\sum_{r=1}^{q}||\mu _{r}(t,\cdot )||^{2}\right] dt  \notag \\
+2p||\tilde{\omega}_{k}(t,\cdot )||^{2(p-1)}\sum_{r=1}^{q}(\mu _{r}(t,\cdot
),\tilde{\omega}_{k}(t,\cdot ))dw_{r}(t),\ t_{k}\wedge \tau _{n}\leq t\leq
t_{k+1}\wedge \tau _{n},  \notag \\
||\tilde{\omega}_{k}(t_{k},\cdot )||^{2p}=E||\tilde{\omega}%
_{k-1}(t_{k},\cdot )||^{2p}.  \notag
\end{gather}%
We have%
\begin{gather*}
dE||\tilde{\omega}_{k}(t,\cdot )||^{2p} \\
=2p\left[ -\frac{\sigma ^{2}}{2}E\left( ||\tilde{\omega}_{k}(t,\cdot
)||^{2(p-1)}||\nabla \tilde{\omega}_{k}(t,\cdot )||^{2}\right) +E||\tilde{%
\omega}_{k}(t,\cdot )||^{2(p-1)}(g(t,\cdot ),\tilde{\omega}_{k}(t,\cdot
))\right. \\
\left. +\frac{2p-1}{2}E||\tilde{\omega}_{k}(t,\cdot
)||^{2(p-1)}\sum_{r=1}^{q}||\mu _{r}(t,\cdot )||^{2}\right] dt,\ t_{k}\wedge
\tau _{n}\leq t\leq t_{k+1}\wedge \tau _{n}, \\
E||\tilde{\omega}_{k}(t_{k},\cdot )||^{2}=E||\tilde{\omega}%
_{k-1}(t_{k},\cdot )||^{2}.
\end{gather*}%
By Poincare's inequality (\ref{Poin}) and doing simple re-arrangements, we
arrive at%
\begin{gather*}
dE||\tilde{\omega}_{k}(t,\cdot )||^{2p} \\
\leq 2p\left[ -\alpha \frac{\sigma ^{2}}{2}E||\tilde{\omega}_{k}(t,\cdot
)||^{2p}+\alpha \frac{\sigma ^{2}}{4}E||\tilde{\omega}_{k}(t,\cdot )||^{2p}+%
\frac{1}{\alpha \sigma ^{2}}||g(t,\cdot )||^{2}\right. \\
\left. +\alpha \frac{\sigma ^{2}}{4}E||\tilde{\omega}_{k}(t,\cdot )||^{2p}+%
\frac{2\left[ (2p-1)(p-1)\right] ^{p}}{\left( \alpha \sigma ^{2}\right)
^{p-1}(p-1)}\left[ \sum_{r=1}^{q}||\mu _{r}(t,\cdot )||^{2}\right] ^{p}%
\right] dt.
\end{gather*}%
We note that the constant $\alpha >0$ in the above expression is due to
Poincare's inequality (\ref{Poin}) and it is, of course, independent of $h$
and $k.$ Hence
\begin{equation}
dE||\tilde{\omega}_{k}(t,\cdot )||^{2p}\leq K\left[ ||g(t,\cdot )||^{2}+%
\left[ \sum_{r=1}^{q}||\mu _{r}(t,\cdot )||^{2}\right] ^{p}\right] dt,
\notag
\end{equation}%
where the constant $K>0$ depends on $p$ but independent of $h$ and $k.$ The
previous inequality implies%
\begin{eqnarray*}
E||\tilde{\omega}_{k}(t_{k+1}\wedge \tau _{n},\cdot )||^{2p} &\leq &E||%
\tilde{\omega}_{k-1}(t_{k},\cdot )||^{2p} \\
&&+KE\int_{t_{k}\wedge \tau _{n}}^{t_{k+1}\wedge \tau _{n}}\left(
||g(s,\cdot )||^{2}+\left[ \sum_{r=1}^{q}||\mu _{r}(s,\cdot )||^{2}\right]
^{p}\right) ds \\
&\leq &E||\phi (\cdot )||^{2p}+E\int_{0}^{t_{k+1}\wedge \tau _{n}}\left(
K||g(s,\cdot )||^{2}+\sum_{r=1}^{q}||\mu _{r}(s,\cdot )||^{2}\right) ds,
\end{eqnarray*}%
and letting $n\rightarrow \infty $ we arrive at (\ref{sns27}). The estimate (%
\ref{sns28}) is evident (see e.g. (\ref{Mnew1})). Theorem~\ref{thm:sns2} is
proved.

\begin{remark}
It is note difficult to see that repeating the proof of Lemma~\ref{LemMat}
word by word, we immediately get that the exponential moment for $||\tilde{%
\omega}_{k}(t_{k+1},\cdot )||^{2}$ is bounded, more precisely the estimate
of the form (\ref{sns30x}) holds for $||\tilde{\omega}_{k}(t_{k+1},\cdot
)||^{2}$ under Assumption~5.1 with $m=0.$
\end{remark}

Now we consider uniform bounds for moments of higher Sobolev norms of $%
\tilde{\omega}_{k}.$

\begin{theorem}
\label{thm:sns2nn}Let Assumption~5.1 hold with $m>0.$ Then
\begin{equation}
E||\tilde{\omega}_{k}(t_{k+1},\cdot )||_{m}^{2p}\leq ||\phi (\cdot
)||_{m}^{2p}+Kt_{k+1},  \label{snd1}
\end{equation}%
where $K>0$ is independent of $h$ and $t_{k}$.
\end{theorem}

\noindent \textbf{Proof}. \ The proof is by induction. To this end, we
assume that moments $E||\tilde{\omega}_{k}(t,\cdot )||_{m-1}^{2p}$ are
bounded (uniformly in $k$ and $h)$ and for sufficiently large $p\geq 1$
(note that Theorem~\ref{thm:sns2} guarantees their boundedness for $m=0).$

We will be adapting recipes from \cite[Section 4.1]{T}. Let the operator $%
\Lambda $ be such that $\Lambda ^{2}=-\Delta .$ We have\ for an integer $%
m\geq 1$ (cf. \cite[p. 29]{T} and also \cite[Section 3.4]{MatPhD}):
\begin{gather}
d||\tilde{\omega}_{k}(t,\cdot )||_{m}^{2p}=2p||\tilde{\omega}_{k}(t,\cdot
)||_{m}^{2(p-1)}\left[ -\frac{\sigma ^{2}}{2}||\tilde{\omega}_{k}(t,\cdot
)||_{m+1}^{2}+(\tilde{\omega}_{k}(t,\cdot ),g(t,\cdot ))_{m}\right.
\label{qqq} \\
\left. -((\hat{v}_{k}(\cdot ),\nabla )\tilde{\omega}_{k}(t,\cdot ),\Lambda
^{2m}\tilde{\omega}_{k}(t,\cdot ))+\frac{2p-1}{2}\sum_{r=1}^{q}||\mu
_{r}(t,\cdot )||_{m}^{2}\right] dt  \notag \\
+2p||\tilde{\omega}_{k}(t,\cdot )||_{m}^{2(p-1)}\sum_{r=1}^{q}(\mu
_{r}(t,\cdot ),\tilde{\omega}_{k}(t,\cdot ))_{m}dw_{r}(t),  \notag \\
\ t_{k}\wedge \tau _{n}\leq t\leq t_{k+1}\wedge \tau _{n},  \notag \\
||\tilde{\omega}_{k}(t_{k},\cdot )||_{m}^{2p}=||\tilde{\omega}%
_{k-1}(t_{k},\cdot )||_{m}^{2p}.  \notag
\end{gather}%
Here $\tau _{n}$ is as in Theorem~\ref{thm:sns2}.

Let us analyze terms in the right-hand side of (\ref{qqq}). We have (e.g.
see \cite[Eq. (4.4)]{T}):%
\begin{eqnarray*}
|(\tilde{\omega}_{k}(t,\cdot ),g(t,\cdot ))_{m}| &\leq &||g(t,\cdot
)||_{m-1}||\tilde{\omega}_{k}(t,\cdot )||_{m+1} \\
&\leq &\frac{4}{\sigma ^{2}}||g(t,\cdot )||_{m-1}^{2}+\frac{\sigma ^{2}}{16}%
||\tilde{\omega}_{k}(t,\cdot )||_{m+1}^{2}
\end{eqnarray*}%
and
\begin{equation*}
K||\tilde{\omega}_{k}(t,\cdot )||_{m}^{2(p-1)}||g(t,\cdot )||_{m-1}^{2}\leq
\frac{K^{p}}{p}\left( \frac{16p}{\alpha \sigma ^{2}(p-1)}\right)
^{p-1}||g(t,\cdot )||_{m-1}^{2p}+\alpha \frac{\sigma ^{2}}{16}||\tilde{\omega%
}_{k}(t,\cdot )||_{m}^{2p},
\end{equation*}%
where as before the constant $\alpha >0$ is due to Poincare's inequality (%
\ref{Poin}). Also, for some $K>0$ dependent on $p:$%
\begin{equation*}
p(2p-1)||\tilde{\omega}_{k}(t,\cdot )||_{m}^{2(p-1)}\sum_{r=1}^{q}||\mu
_{r}(t,\cdot )||_{m}^{2}\leq K\left( \sum_{r=1}^{q}||\mu _{r}(t,\cdot
)||_{m}^{2}\right) ^{p}+\alpha \frac{\sigma ^{2}}{8}||\tilde{\omega}%
_{k}(t,\cdot )||_{m}^{2p}.
\end{equation*}%
Hence, we can write for some $K>0:$%
\begin{gather*}
d||\tilde{\omega}_{k}(t,\cdot )||_{m}^{2p}\leq \left[ -\frac{\sigma ^{2}}{4}%
p\left\{ ||\tilde{\omega}_{k}(t,\cdot )||_{m}^{2(p-1)}||\tilde{\omega}%
_{k}(t,\cdot )||_{m+1}^{2}+\alpha ||\tilde{\omega}_{k}(t,\cdot
)||_{m}^{2p}\right\} \right.  \\
+K||g(t,\cdot )||_{m-1}^{2p}-2p\left\{ ||\tilde{\omega}_{k}(t,\cdot
)||_{m}^{2(p-1)}((\hat{v}_{k}(\cdot ),\nabla )\tilde{\omega}_{k}(t,\cdot
),\Lambda ^{2m}\tilde{\omega}_{k}(t,\cdot ))\right\}  \\
\left. +K\left( \sum_{r=1}^{q}||\mu _{r}(t,\cdot )||_{m}^{2}\right) ^{p}
\right] dt \\
+2p||\tilde{\omega}_{k}(t,\cdot )||_{m}^{2(p-1)}\sum_{r=1}^{q}(\mu
_{r}(t,\cdot ),\tilde{\omega}_{k}(t,\cdot ))_{m}dw_{r}(t) \\
\ t_{k}\wedge \tau _{n}\leq t\leq t_{k+1}\wedge \tau _{n},\ \ ||\tilde{\omega%
}_{k}(t_{k},\cdot )||_{m}^{2p}=||\tilde{\omega}_{k-1}(t_{k},\cdot
)||_{m}^{2p}.
\end{gather*}%
Let us now estimate the trilinear-form:
\begin{gather*}
|((\hat{v}_{k}(\cdot ),\nabla )\tilde{\omega}_{k}(t,\cdot ),\Lambda ^{2m}%
\tilde{\omega}_{k}(t,\cdot ))|\leq K\sum_{l=1}^{m}||\hat{v}_{k}(\cdot
)||_{l}||\tilde{\omega}_{k}(t,\cdot )||_{m-l+3/2}||\tilde{\omega}%
_{k}(t,\cdot )||_{m+1} \\
\leq K||\tilde{\omega}_{k}(t_{k},\cdot )||_{m-1}||\tilde{\omega}_{k}(t,\cdot
)||_{m+1/2}||\tilde{\omega}_{k}(t,\cdot )||_{m+1} \\
\leq K||\tilde{\omega}_{k}(t,\cdot )||_{m-1}^{5/4}||\tilde{\omega}%
_{k}(t,\cdot )||_{m+1}^{7/4} \\
\leq \frac{K^{8}}{8}\left( \frac{56}{\sigma ^{2}}\right) ^{7}||\tilde{\omega}%
_{k}(t_{k},\cdot )||_{m-1}^{10}+\frac{\sigma ^{2}}{8}||\tilde{\omega}%
_{k}(t,\cdot )||_{m+1}^{2},
\end{gather*}%
where for the first line we used the recipe from \cite[pp. 29-30]{T} and the
inequality (\ref{A2}); for the second line we used (\ref{A3}) and that $%
||u(\cdot )||_{m_{1}}\leq ||u(\cdot )||_{m_{2}}$ for $m_{2}\geq m_{1}$; the
third line is obtained using (\ref{A4i}); and the fourth line follows from
Young's inequality. Note that the constants $K>0$ in the first and second
lines are different.

Further,
\begin{eqnarray*}
&&K||\tilde{\omega}_{k}(t,\cdot )||_{m-1}^{10}||\tilde{\omega}_{k}(t,\cdot
)||_{m}^{2(p-1)} \\
&\leq &\frac{K^{p}}{p}\left( \frac{8p}{\alpha \sigma ^{2}(p-1)}\right)
^{p-1}||\tilde{\omega}_{k}(t_{k},\cdot )||_{m-1}^{10p}+\alpha \frac{\sigma
^{2}}{8}||\tilde{\omega}_{k}(t,\cdot )||_{m}^{2p}.
\end{eqnarray*}%
Thus, for some $K>0$%
\begin{gather*}
dE||\tilde{\omega}_{k}(t,\cdot )||_{m}^{2p}\leq  \\
+K\left[ ||g(t,\cdot )||_{m-1}^{2p}+E||\tilde{\omega}_{k}(t_{k},\cdot
)||_{m-1}^{10p}+\left( \sum_{r=1}^{q}||\mu _{r}(t,\cdot )||_{m}^{2}\right)
^{p}\right] dt, \\
\ t_{k}\wedge \tau _{n}\leq t\leq t_{k+1}\wedge \tau _{n},\ \ E||\tilde{%
\omega}_{k}(t_{k},\cdot )||_{m}^{2p}=E||\tilde{\omega}_{k-1}(t_{k},\cdot
)||_{m}^{2p}.
\end{gather*}%
By the Cauchy-Bunyakovskii inequality and the induction assumption at the
start of the proof, we get%
\begin{equation*}
E||\tilde{\omega}_{k}(t_{k},\cdot )||_{m-1}^{10p}\leq K
\end{equation*}%
with a constant $K>0$ independent of $h$ and $k.$ Hence
\begin{gather*}
dE||\tilde{\omega}_{k}(t,\cdot )||_{m}^{2p}\leq Kdt, \\
\ t_{k}\wedge \tau _{n}\leq t\leq t_{k+1}\wedge \tau _{n},\ \ E||\tilde{%
\omega}_{k}(t_{k},\cdot )||_{m}^{2p}=E||\tilde{\omega}_{k-1}(t_{k},\cdot
)||_{m}^{2p}
\end{gather*}%
and
\begin{equation*}
E||\tilde{\omega}_{k}(t_{k+1}\wedge \tau _{n},\cdot )||_{m}^{2p}\leq E||%
\tilde{\omega}_{k-1}(t_{k},\cdot )||_{m}^{2p}+Kh,
\end{equation*}%
from which (\ref{snd1}) follows by the standard arguments. Theorem~\ref%
{thm:sns2nn} is proved.

\subsection{Mean-square convergence theorem\label{sec:sns4}}

To prove the global convergence of $\tilde{\omega}_{k}(t_{k},\cdot ),$ we
use the idea of the proof of the fundamental theorem of mean-square
convergence for SDEs \cite{Mil87} (see also \cite[Section 1.1]{MT1}).

\begin{theorem}
\label{thm:sns3}Let Assumption~5.1 hold with $m=2.$ The global approximation
$\tilde{\omega}_{k}(t_{k+1},x)$ for the problem $(\ref{sns7})$-$(\ref{sns9})$
has the first mean-square order accuracy.
\end{theorem}

\noindent \textbf{Proof}. \ We note that in the proof we shall again use
letters $K$ and $C(\cdot ,\mathbf{\omega })$ to denote various deterministic
constants and random variables, respectively, which are independent of $h$
and $k,$ and $K$ is also independent of $h$ and $k$; their values may change
from line to line.

We have
\begin{eqnarray}
R(t_{k+1},x) &:&=\omega (t_{k+1},x;0,\phi )-\tilde{\omega}%
_{k}(t_{k+1},x;0,\phi )  \label{Ba27} \\
&=&\omega (t_{k+1},x;t_{k},\omega (t_{k},\cdot ))-\tilde{\omega}%
_{k}(t_{k+1},x;t_{k},\tilde{\omega}_{k}(t_{k},\cdot ))  \notag \\
&=&(\omega (t_{k+1},x;t_{k},\omega (t_{k},\cdot ))-\omega (t_{k+1},x;t_{k},%
\tilde{\omega}_{k}(t_{k},\cdot )))  \notag \\
&&+(\omega (t_{k+1},x;t_{k},\tilde{\omega}_{k}(t_{k},\cdot ))-\tilde{\omega}%
_{k}(t_{k+1},x;t_{k},\tilde{\omega}_{k}(t_{k},\cdot )))\,,  \notag
\end{eqnarray}%
where $\cdot $ reflects function dependence of solutions on the initial
conditions. The first difference in the right-hand side of (\ref{Ba27}) is
the error of the solution arising due to the error in the initial data at
time $t_{k},$ accumulated at the $k$-th step. The second difference is the
one-step error at the $(k+1)$-step:
\begin{equation}
\delta _{\omega }(t_{k+1},x):=\omega (t_{k+1},x;t_{k},\tilde{\omega}%
_{k}(t_{k},\cdot ))-\tilde{\omega}_{k}(t_{k+1},x;t_{k},\tilde{\omega}%
_{k}(t_{k},\cdot ))  \label{Ba30}
\end{equation}%
for which estimates are given in Theorem~\ref{thm:sns1} taking into account
that Theorems~\ref{thm:sns2} and \ref{thm:sns2nn} guarantees boundedness of
moments of $||\tilde{\omega}_{k}(t_{k},\cdot )||_{3}$ under the conditions
of this theorem. Taking the $L^{2}$-norm of both sides of (\ref{Ba27}), we
obtain
\begin{gather}
||R(t_{k+1},\cdot )||^{2}=||\omega (t_{k+1},\cdot ;t_{k},\omega (t_{k},\cdot
))-\omega (t_{k+1},\cdot ;t_{k},\tilde{\omega}_{k}(t_{k},\cdot ))||^{2}\
\label{Ba28} \\
+||\delta _{\omega }(t_{k+1},\cdot )||^{2}+2(\omega (t_{k+1},\cdot
;t_{k},\omega (t_{k},\cdot ))-\omega (t_{k+1},\cdot ;t_{k},\tilde{\omega}%
_{k}(t_{k},\cdot )),\delta _{\omega }(t_{k+1},\cdot )),  \notag
\end{gather}%
where the first $\cdot $ in each $\omega $ or $\tilde{\omega}_{k}$ reflects
that we took $L^{2}$-norm.

Using (\ref{sns19}) from Theorem~\ref{thm:sns1} together with Theorems~\ref%
{thm:sns2} and \ref{thm:sns2nn}, we obtain for the second term in (\ref{Ba28}%
):
\begin{equation}
||\delta _{\omega }(t_{k+1},\cdot )||^{2}\leq C(t_{k+1},\mathbf{\omega }%
)h^{3},  \label{Ba288}
\end{equation}%
where $C(t_{k+1},\mathbf{\omega })>0$ is an $\mathcal{F}_{t_{k+1}}$%
-measurable with bounded second moment.

By (\ref{sns31}) from Lemma \ref{Lem2} together with Theorems~\ref{thm:sns2}
and \ref{thm:sns2nn}, we get for the first term in (\ref{Ba28}):%
\begin{eqnarray}
&&||\omega (t_{k+1},x;t_{k},\omega (t_{k},\cdot ))-\omega (t_{k+1},x;t_{k},%
\tilde{\omega}_{k}(t_{k},\cdot ))||^{2}  \label{Ba29} \\
&\leq &||R(t_{k},\cdot )||^{2}\exp \left( Kh+c\int_{t_{k}}^{t_{k+1}}||\nabla
\omega (s^{\prime },\cdot ;t_{k},\omega (t_{k},\cdot ))||^{2}ds^{\prime
}\right) .  \notag
\end{eqnarray}%
The difference $\omega (t_{k+1},x;t_{k},\omega (t_{k},\cdot ))-\omega
(t_{k+1},x;t_{k},\tilde{\omega}_{k}(t_{k},\cdot ))$ in the last summand in (%
\ref{Ba28}) can be treated using (\ref{sns30}) from Lemma~\ref{Lem2}:
\begin{equation*}
\omega (t_{k+1},x;t_{k},\omega (t_{k},\cdot ))-\omega (t_{k+1},x;t_{k},%
\tilde{\omega}_{k}(t_{k},\cdot ))=R(t_{k},x)+\eta (t_{k},x)\,.
\end{equation*}%
Using a conditional version of (\ref{sns19}) from Theorem~\ref{thm:sns1} and
(\ref{sns32}) from Lemma~\ref{Lem2} together with Theorems~\ref{thm:sns2}
and \ref{thm:sns2nn}, we get
\begin{eqnarray}
|((\eta (t_{k},\cdot ),\delta _{\omega }(t_{k+1},\cdot ))| &\leq &||\eta
(t_{k},\cdot )||||\delta _{\omega }(t_{k+1},\cdot )||\leq ||\eta
(t_{k},\cdot )||^{2}+\frac{1}{4}||\delta _{\omega }(t_{k+1},\cdot )||^{2}\ \
\ \ \   \label{Ba32} \\
&\leq &h||R(t_{k},\cdot )||^{2}+C(t_{k+1},\mathbf{\omega })h^{3},  \notag
\end{eqnarray}%
where $C(t_{k+1},\mathbf{\omega })>0$ is an $\mathcal{F}_{t_{k+1}}$%
-measurable with bounded second moment.

Combining the above, we arrive at
\begin{eqnarray}
||R(t_{k+1},\cdot )||^{2} &\leq &|R(t_{k},\cdot )||^{2}\exp \left(
Kh+c\int_{t_{k}}^{t_{k+1}}||\nabla \omega (s^{\prime },\cdot ;t_{k},\omega
(t_{k},\cdot ))||^{2}ds^{\prime }\right)  \label{bax1} \\
&&+h||R(t_{k},\cdot )||^{2}+(R(t_{k},\cdot ),\delta _{\omega }(t_{k+1},\cdot
))+C(t_{k+1},\mathbf{\omega })h^{3}.  \notag
\end{eqnarray}

Since $||R(0,\cdot )||=0,$ summing (\ref{bax1}) from $k=0$ to $n$, we get
\begin{eqnarray*}
&&||R(t_{n+1},\cdot )||^{2} \\
&\leq &\sum_{k=1}^{n}||R(t_{k},\cdot )||^{2}\left[ \exp \left(
Kh+c\int_{t_{k}}^{t_{k+1}}||\nabla \omega (s^{\prime },\cdot ;t_{k},\omega
(t_{k},\cdot ))||^{2}ds^{\prime }\right) -1+h\right]  \\
&&+h^{3}\sum_{k=0}^{n}C(t_{k+1},\mathbf{\omega })+\sum_{k=1}^{n}(R(t_{k},%
\cdot ),\delta _{\omega }(t_{k+1},\cdot )) \\
&\leq &\sum_{k=1}^{n}||R(t_{k},\cdot )||^{2}\left[ \exp \left(
Kh+c\int_{t_{k}}^{t_{k+1}}||\nabla \omega (s^{\prime },\cdot ;t_{k},\omega
(t_{k},\cdot ))||^{2}ds^{\prime }\right) -1\right]  \\
&&+h^{3}\sum_{k=0}^{n}C(t_{k+1},\mathbf{\omega })+\sum_{k=1}^{n}(R(t_{k},%
\cdot ),\delta _{\omega }(t_{k+1},\cdot )).
\end{eqnarray*}%
From which, by a version of Gronwall's lemma (see, e.g. \cite{Gronw,Kruse18}%
), we obtain
\begin{eqnarray}
||R(t_{n+1},\cdot )||^{2} &\leq &F_{n}+\sum_{k=1}^{n}F_{k-1}  \label{bax3} \\
&&\cdot \left[ \exp \left( Kh+c\int_{t_{k}}^{t_{k+1}}||\nabla \omega
(s^{\prime },\cdot ;t_{k},\omega (t_{k},\cdot ))||^{2}ds^{\prime }\right) -1%
\right]   \notag \\
&&\cdot \mathop{\displaystyle \prod }\limits_{j=k+1}^{n}\exp \left(
Kh+c\int_{t_{j}}^{t_{j+1}}||\nabla \omega (s^{\prime },\cdot ;t_{k},\omega
(t_{k},\cdot ))||^{2}ds^{\prime }\right) ,  \notag
\end{eqnarray}%
where
\begin{equation*}
F_{k}:=h^{3}\sum_{j=0}^{k}C(t_{j+1},\mathbf{\omega })+%
\sum_{j=1}^{k}(R(t_{j},\cdot ),\delta _{\omega }(t_{j+1},\cdot )).
\end{equation*}%
We have
\begin{eqnarray}
&&||R(t_{n+1},\cdot )||^{2}\leq F_{n}  \label{bax31} \\
&&+\sum_{k=1}^{n}F_{k-1}\cdot \left[ \exp \left(
Kh+c\int_{t_{k}}^{t_{k+1}}||\nabla \omega (s^{\prime },\cdot ;t_{k},\omega
(t_{k},\cdot ))||^{2}ds^{\prime }\right) -1\right]   \notag \\
&&\cdot \exp \left( K(t_{n+1}-t_{k+1})+c\int_{t_{k+1}}^{t_{n+1}}||\nabla
\omega (s^{\prime },\cdot ;t_{k+1},\omega (t_{k+1},\cdot ))||^{2}ds^{\prime
}\right)   \notag \\
&=&F_{n}+\sum_{k=1}^{n}F_{k-1}\cdot \left[ \exp \left(
K(t_{n+1}-t_{k})+c\int_{t_{k}}^{t_{n+1}}||\nabla \omega (s^{\prime },\cdot
;t_{k},\omega (t_{k},\cdot ))||^{2}ds^{\prime }\right) \right.   \notag \\
&&\left. -\exp \left( K(t_{n+1}-t_{k+1})+c\int_{t_{k+1}}^{t_{n+1}}||\nabla
\omega (s^{\prime },\cdot ;t_{k+1},\omega (t_{k+1},\cdot ))||^{2}ds^{\prime
}\right) \right]   \notag \\
&=&\sum_{k=1}^{n}\left( F_{k}-F_{k-1}\right) \exp \left(
K(t_{n+1}-t_{k+1})+c\int_{t_{k+1}}^{t_{n+1}}||\nabla \omega (s^{\prime
},\cdot ;t_{k+1},\omega (t_{k+1},\cdot ))||^{2}ds^{\prime }\right)   \notag
\\
&&+h^{3}C(t_{1},\mathbf{\omega })\exp \left(
K(t_{n+1}-t_{1})+c\int_{t_{1}}^{t_{n+1}}||\nabla \omega (s^{\prime },\cdot
;t_{1},\omega (t_{1},\cdot ))||^{2}ds^{\prime }\right) .  \notag
\end{eqnarray}

For the last term in the right-hand side of (\ref{bax31}), we obtain using
the Cauchy-Bunyakovsky inequality and Lemma~\ref{LemMat}:%
\begin{equation}
E\left\{ h^{3}C(t_{1},\mathbf{\omega })\exp \left(
K(t_{n+1}-t_{1})+c\int_{t_{k+1}}^{t_{n+1}}||\nabla \omega (s^{\prime },\cdot
;0,\phi (\cdot ))||^{2}ds^{\prime }\right) \right\} \leq Kh^{3}.
\label{bax5}
\end{equation}

Consider now the first term in the right-hand side of (\ref{bax31}). We have
\begin{eqnarray}
&&\left( F_{k}-F_{k-1}\right) \exp \left(
K(t_{n+1}-t_{k+1})+c\int_{t_{k+1}}^{t_{n+1}}||\nabla \omega (s^{\prime
},\cdot ;t_{k+1},\omega (t_{k+1},\cdot ))||^{2}ds^{\prime }\right) \ \ \ \ \
\label{bax55} \\
&=&\exp \left( K(t_{n+1}-t_{k+1})+c\int_{t_{k+1}}^{t_{n+1}}||\nabla \omega
(s^{\prime },\cdot ;t_{k+1},\omega (t_{k+1},\cdot ))||^{2}ds^{\prime
}\right)   \notag \\
&&\times \left[ h^{3}C(t_{k+1},\mathbf{\omega })+(R(t_{k},\cdot ),\delta
_{\omega }(t_{k+1},\cdot ))\right] .  \notag
\end{eqnarray}%
Expectation of the first term from the right-hand side of (\ref{bax55}) is
estimated by $Kh^{3}$ as in (\ref{bax5}). Let us now consider the second
term.

By the martingale representation theorem and Lemma~\ref{LemMat}, we can
obtain
\begin{eqnarray}
&&E\left[ \left. \exp \left(
K(t_{n+1}-t_{k+1})+c\int_{t_{k+1}}^{t_{n+1}}||\nabla \omega (s^{\prime
},\cdot ;t_{k+1},\omega (t_{k+1},\cdot ))||^{2}ds^{\prime }\right)
\right\vert \mathcal{F}_{t_{k+1}}\right] \ \ \ \ \ \   \label{snsmrt} \\
&=&E\left[ \exp \left( K(t_{n+1}-t_{k+1})+c\int_{t_{k+1}}^{t_{n+1}}||\nabla
\omega (s^{\prime },\cdot ;t_{k+1},\omega (t_{k+1},\cdot ))||^{2}ds^{\prime
}\right) \right]  \notag \\
&&+\sum_{r=1}^{q}\int_{0}^{t_{k+1}}\lambda _{r}(s)dw_{r}(s),  \notag
\end{eqnarray}%
where $\lambda _{r}(s)$ are $\mathcal{F}_{s}$-adapted square-integrable
stochastic processes.

Using (\ref{sns18}) and a conditional version of (\ref{sns19}) from Theorem~%
\ref{thm:sns1} together with Theorems~\ref{thm:sns2} and \ref{thm:sns2nn}
and also using Lemma~\ref{LemMat} and (\ref{snsmrt}), we arrive at
\begin{eqnarray*}
&&\left\vert E\left\{ \exp \left(
K(t_{n+1}-t_{k+1})+c\int_{t_{k+1}}^{t_{n+1}}||\nabla \omega (s^{\prime
},\cdot ;t_{k+1},\omega (t_{k+1},\cdot ))||^{2}ds^{\prime }\right) \right.
\right.  \\
&&\left. \left. \Bigl.\overset{\ }{\ \ }(R(t_{k},\cdot ),\delta _{\omega
}(t_{k+1},\cdot ))\Bigr\vert\mathcal{F}_{t_{k}}\right\} \right\vert  \\
&=&\left\vert \left( R(t_{k},\cdot ),E\left\{ \exp \left(
K(t_{n+1}-t_{k+1})+c\int_{t_{k+1}}^{t_{n+1}}||\nabla \omega (s^{\prime
},\cdot ;t_{k+1},\omega (t_{k+1},\cdot ))||^{2}ds^{\prime }\right) \right.
\right. \right.  \\
&&\left. \left. \left. \Biggl.\overset{\ }{\ \ }\delta _{\omega
}(t_{k+1},\cdot )\Biggr\vert\mathcal{F}_{t_{k}}\right\} \right) \right\vert
\\
&\leq &||R(t_{k},\cdot )||\left\Vert E\left\{ \exp \left(
K(t_{n+1}-t_{k+1})+c\int_{t_{k+1}}^{t_{n+1}}||\nabla \omega (s^{\prime
},\cdot ;t_{k+1},\omega (t_{k+1},\cdot ))||^{2}ds^{\prime }\right) \right.
\right.  \\
&&\left. \left. \left. \Biggl.\delta _{\omega }(t_{k+1},\cdot )\Biggr\vert%
\mathcal{F}_{t_{k}}\right\} \right) \right\Vert
\end{eqnarray*}%
\begin{eqnarray*}
&=&||R(t_{k},\cdot )|| \\
&&\cdot \left\Vert E\left\{ \left( E\left[ \exp \left(
K(t_{n+1}-t_{k+1})+c\int_{t_{k+1}}^{t_{n+1}}||\nabla \omega (s^{\prime
},\cdot ;t_{k+1},\omega (t_{k+1},\cdot ))||^{2}ds^{\prime }\right) \right]
\right. \right. \right.  \\
&&\left. \left. \left. \left. +\sum_{r=1}^{q}\int_{0}^{t_{k+1}}\lambda
_{r}(s)dw_{r}\right) \delta _{\omega }(t_{k+1},\cdot )\right\vert \mathcal{F}%
_{t_{k}}\right\} \right\Vert  \\
&\leq &||R(t_{k},\cdot )||E\left[ \exp \left(
K(t_{n+1}-t_{k+1})+c\int_{t_{k+1}}^{t_{n+1}}||\nabla \omega (s^{\prime
},\cdot ;t_{k+1},\omega (t_{k+1},\cdot ))||^{2}ds^{\prime }\right) \right]
\\
&&\cdot \left\Vert E\left\{ \left. \delta _{\omega }(t_{k+1},\cdot
)\right\vert \mathcal{F}_{t_{k}}\right\} \right\Vert +||R(t_{k},\cdot
)||\left\Vert E\left\{ \left. \delta _{\omega }(t_{k+1},\cdot
)\sum_{r=1}^{q}\int_{0}^{t_{k+1}}\lambda _{r}(s)dw_{r}\right\vert \mathcal{F}%
_{t_{k}}\right\} \right\Vert
\end{eqnarray*}%
\begin{eqnarray*}
&\leq &||R(t_{k},\cdot )||C(t_{k},\mathbf{\omega })h^{2}+||R(t_{k},\cdot
)||\left\Vert \sum_{r=1}^{q}\int_{0}^{t_{k}}\lambda _{r}(s)dw_{r}E\left\{
\left. \delta _{\omega }(t_{k+1},\cdot )\right\vert \mathcal{F}%
_{t_{k}}\right\} \right\Vert  \\
&&+||R(t_{k},\cdot )||\left( E\left[ ||\left. \delta _{\omega
}(t_{k+1},\cdot )||^{2}\right\vert \mathcal{F}_{t_{k}}\right] \right) ^{1/2}
\\
&&\cdot \left( \sum_{r=1}^{q}\left( E\left[ \left. \left\{
\int_{t_{k}}^{t_{k+1}}\lambda _{r}(s)dw_{r}\right\} ^{2}\right\vert \mathcal{%
F}_{t_{k}}\right] \right) ^{1/2}\right)  \\
&\leq &||R(t_{k},\cdot )||C(t_{k},\mathbf{\omega })h^{2}\leq
h||R(t_{k},\cdot )||^{2}+\frac{h^{3}}{4}C^{2}(t_{k},\mathbf{\omega }).
\end{eqnarray*}%
Therefore,
\begin{eqnarray}
&&E\left\{ \left( F_{k}-F_{k-1}\right) \exp \left(
K(t_{n+1}-t_{k+1})+c\int_{t_{k+1}}^{t_{n+1}}||\nabla \omega (s^{\prime
},\cdot ;t_{k+1},\omega (t_{k+1},\cdot ))||^{2}ds^{\prime }\right) \right\}
\ \ \ \ \ \ \ \ \ \ \   \label{bax99} \\
&\leq &hE||R(t_{k},\cdot )||^{2}+Kh^{3}.  \notag
\end{eqnarray}%
From (\ref{bax31}), (\ref{bax5}), and (\ref{bax99}), we obtain
\begin{equation*}
E||R(t_{n+1},\cdot )||^{2}\leq h\sum_{k=1}^{n}E||R(t_{k},\cdot
)||^{2}+Kh^{2},
\end{equation*}%
from which it follows by a version of Gronwall's lemma that
\begin{equation*}
E||R(t_{n+1},\cdot )||^{2}\leq Kh^{2}
\end{equation*}%
as required. Theorem~\ref{thm:sns3} is proved.

\begin{remark}
Various approaches can be used to turn the method, introduced at the start
of Section~\ref{sec:sns3} for the problem $(\ref{sns7})$-$(\ref{sns9})$,
into a numerical algorithm. To obtain a constructive numerical algorithm, we
need to approximate the linear SPDE $(\ref{sns11})$-$(\ref{sns12})$ at every
step. To this end, for instance, we can discretize this SPDE in space using
the spectral method based on the Fourier expansion and use a finite
difference for time discretization (see such an algorithm in the
deterministic setting in e.g. \cite{Peyret}). Alternatively, we can apply
the method based on averaging characteristics to $(\ref{sns11})$-$(\ref%
{sns12})$ \cite{spde}. We leave construction, analysis and testing of such
algorithms for a future work.
\end{remark}


\begin{thebibliography}{99}
\bibitem{Gronw} P.R. Beesack. More generalised discrete Gronwall
inequalities. \textit{ZAMM Z. Angew. Math. Mech.} \textbf{65} (1985),
589--595.

\bibitem{BBM14} H. Bessaih, Z. Brze\'{z}niak, A. Millet. Splitting up method
for the 2D stochastic Navier--Stokes equations. \textit{SPDEs: An. Comp. }%
\textbf{2} (2014), 433--470.

\bibitem{BCP10} Z. Brze\'{z}niak, E. Carelli, A. Prohl. Finite element based
discretizations of the incompressible Navier-Stokes equations with
multiplicative random forcing. \textit{IMA J. Num. Anal.} \textbf{34}
(2014), 502--549.

\bibitem{BFR} B. Busnello, F. Flandoli, M. Romito. A probabilistic
representation for the vorticity of a $3D$ viscous fluid and for general
systems of parabolic equations. \textit{Proc. Edinb. Math. Soc.} \textbf{48 }%
(2005), 295--336.

\bibitem{CP} E. Carelli, A. Prohl. Rates of convergence for discretizations
of the stochastic incompressible Navier-Stokes equations. \textit{SIAM J.
Num. Anal.} \textbf{50} (2012), 2467--2496.

\bibitem{CF} P. Constantin, C. Foias. \textit{Navier-Stokes Equations.}
University of Chicago Press, 1988.

\bibitem{Dorsek} P. D\"{o}rsek. Semigroup splitting and cubature
approximations for the stochastic Navier-Stokes equations. \textit{SIAM J.
Num. Anal.} \textbf{50} (2012), 729-746.

\bibitem{DJT} T. Dubois, F. Jauberteau, R. Temam. \textit{Dynamic Multilevel
Methods and the Numerical Simulation of Turbulence}. Cambridge Univ. Press,
1999.

\bibitem{Flandoli} F. Flandoli. An introduction to 3D stochastic fluid
dynamics. In: SPDE in Hydrodynamic: Recent Progress and Prospects, \textit{%
Lecture Notes in Mathematics} \textbf{1942}, Springer, 2008, 51--150.

\bibitem{FMRT01} C. Foias, O. Manley, R. Rosa, R. Temam. \textit{%
Navier-Stokes Equations and Turbulence}. Cambridge University Press, 2001.

\bibitem{FEM} V. Girault, P.A. Raviart. \textit{Finite Element Methods for
Navier-Stokes Equations.} Springer, 1986.

\bibitem{HaMa06} M. Hairer, J.C. Mattingly. Ergodicity of the 2D
Navier-Stokes equations with degenerate stochastic forcing. \textit{Ann.
Math.} \textbf{164} (2006), 993--1032.

\bibitem{Kruse18} R. Kruse, M. Scheutzow. A discrete stochastic Gronwall
lemma. \textit{Math. Comp. Simul.} \textbf{143} (2018), 149--157.

\bibitem{LM} J.L. Lions,\ E. Magenes. \textit{Nonhomogeneous Boundary Value
Problems and Applications}. Springer, 1972.

\bibitem{MB} A.J. Majda, A.L. Bertozzi. \textit{Vorticity and Incompressible
Flow}. Cambridge Univ. Press, 2003.

\bibitem{MatPhD} J.C. Mattingly. \textit{The Stochastic Navier-Stokes
Equation: Energy Estimates and Phase Space Contraction}. Ph.D. Thesis,
Princeton University, Princeton, 1998.

\bibitem{Mat02} J.C. Mattingly. The dissipative scale of the stochastics
Navier--Stokes equation: regularization and analyticity. \textit{J. Stat.
Phys.} \textbf{108} (2002), 1157--1179.

\bibitem{RozNS05} R. Mikulevicius, B. Rozovskii. Global L2-solutions of
stochastic Navier-Stokes equations. \textit{Ann. Prob.} \textbf{33} (2005),
137--176.

\bibitem{M} G.N. Milstein. Probabilistic solution of linear systems of
elliptic and parabolic equations. Theor. Prob. Appl. \textbf{23 }(1978),
851--855.

\bibitem{Mil87} G.N. Milstein. A theorem on the order of convergence of
mean-square approximations of solutions of systems of stochastic
differential equations. \textit{Teor. Veroyat. Primenen}., \textbf{32}
(1987), 809--811.

\bibitem{M1} G.N. Milstein. The probability approach to numerical solution
of nonlinear parabolic equations. \textit{Num. Meth. PDE} \textbf{18}
(2002), 490--522.

\bibitem{MSS} G.N. Milstein, J.G.M. Schoenmakers, V. Spokoiny. Transition
density estimation for stochastic differential equations via forward-reverse
representations. \textit{Bernoulli} \textbf{10 }(2004), 281--312.

\bibitem{MT1} G.N. Milstein, M.V. Tretyakov. \textit{Stochastic Numerics for
Mathematical Physics}. Springer, 2004.

\bibitem{spde} G.N. Milstein, M.V. Tretyakov. Solving parabolic stochastic
partial differential equations via averaging over characteristics. \textit{%
Math. Comp.} \textbf{78} (2009), 2075--2106.

\bibitem{MTns13} G.N. Milstein, M.V. Tretyakov. Layer methods for the
incompressible Navier-Stokes equations with space periodic conditions.
\textit{Adv. App. Prob.} \textbf{45} (2013), 742--772.

\bibitem{MTsns} G.N. Milstein, M.V. Tretyakov. Layer methods for
Navier-Stokes equations with additive noise using simplest characteristics.
\textit{J. Comp Appl. Maths.} \textbf{302} (2016), 1--23.

\bibitem{Peyret} R. Peyret. \textit{Spectral Methods for Incompressible
Viscous Flow}. Springer, 2002.

\bibitem{RT} R. Temam. \textit{Navier-Stokes Equations,\ Theory and
Numerical Analysis.} AMS Chelsea Publishing, 2001.

\bibitem{T} R. Temam. \textit{Navier-Stokes Equations and Nonlinear
Functional Analysis}. SIAM, 1995.

\bibitem{W} P. Wesseling. \textit{Principles of Computational Fluid Dynamics}%
. Springer, 2001.
\end{thebibliography}
\end{document}